\documentclass[12pt]{article}
\usepackage{amsmath,amssymb}

\newtheorem{theorem}{Theorem}
\newtheorem{lemma}{Lemma}
\newtheorem{corollary}{Corollary}
\newtheorem{proposition}{Proposition} 
\newtheorem{defin}{Definition}
\newtheorem{rem}{Remark}

\newcommand{\Bx}{[B^-\times B^+]}
\newcommand{\Di}{\mathfrak D}
\newcommand{\eb}{\hfill q.e.d.}
\newcommand{\borel}{\mathfrak b}

\newcommand{\mor}{morphism}

\newcommand{\x}{\times}
\newcommand{\<}{\langle}
\renewcommand{\>}{\rangle}
\renewcommand{\a}{\alpha}
\renewcommand{\b}{\mathfrak b}
\newcommand{\be}{\beta}
\newcommand{\bs}{\bigskip}
\newcommand{\dt}{\cdot}
\newcommand{\C}{\mathbb C}

\newcommand{\D}{\Delta}
\newcommand{\e}{\varepsilon}
\newcommand{\f}{\mathfrak h}
\newcommand{\g}{\mathfrak g}
\newcommand{\Ge}{\tilde{g}}
\newcommand{\h}{\tilde h}
\newcommand{\G}{\mathfrak g}

\newcommand{\mb}{\mbox}
\newcommand{\N}{\mathfrak n}

\renewcommand{\O}{\omega}
\renewcommand{\o}{\otimes}
\newcommand{\p}{\partial}

\newcommand{\sig}{\sigma}

\newcommand{{\z}}{\Bbb Z}
\newcommand{\M}{\mathcal M}
\newcommand{\xii}{\tilde\xi}
\newcommand{\stk}{\stackrel}

\begin{document}

\title{Integrability of characteristic Hamiltonian systems on simple
Lie groups with standard Poisson Lie structure}

\author{ N.Reshetikhin}
\maketitle

\begin{abstract}
Phase space of a characteristic Hamiltonian system is a symplectic
leaf of a factorizable Poisson Lie group. Its Hamiltonian is a
restriction to the symplectic leaf of a function on the group
which is invariant with respect to conjugations. It is shown in
this paper that such system is always integrable.
\end{abstract}

\tableofcontents
\section{Introduction}

The discovery of Lax pair for the Korteveg-de-Vris equation opened
new class of exactly solvable ordinary and partial differential
equations. Solutions to such equations can be expressed in terms
of solutions to certain spectral problems, or in terms of certain
factorization problems of either Gauss type or Riemann-Hilbert
type.

It has been discovered later that for most of such systems there
is a Hamiltonian structure with respect to which they are
integrable system in a sense of Lioville.

A Lie theoretical explanation of this fact has been provided by
Kostant \cite{Kostant-oxford} on the example of the Toda system.
He noticed that the phase space of the Toda system can be
naturally identified with an orbit of a Borel subgroup acting by
adjoint representation a simple Lie group which passes through the
proncipal nilptent element in the opposite Borel subalgebra. The
Hamiltonians of the Toda system are restrictions of central
functions on the corresponding simple Lie algebra to these orbits.
Then Adler \cite{Adler-1979} applied this approach to the KdV
equation and interpreted it in terms of Lie algebra of
pseudo-differential operators (see also \cite{Simes}). This
approach has been generalized in \cite{R-STS}. The key observation
was that one should consider certain Lie algebra structure on the
space dual to the Lie algebra instead of Borel subalgebras.

It has been discovered in a number of examples \cite{S} that if
the Poisson brackets between matrix elements of Lax operator have
certain special structure (so-called $r$-matrix structure ) and
that spectral functions on such Lax operators generate integrable
systems. In contrast to Toda systems, KdV equations and others
interpreted in terms of Lie algebras, most of these systems had
nonlinear (quadratic) Poisson brackets. Such Poisson structures
were later categorized by Drinfeld \cite{Drinfeld87} when he
introduced the notion of a Poisson Lie group. The theory of
Kostant was generalized to Poisson Lie groups by
Semenov-Tian-Shanski \cite{STS}. In this case the phase space of a
system is a symplectic leaf of a factorizable Poisson Lie group.
Integrals of such system are adjoint invariant functions
restricted to the symplectic leaf. They Poisson commute, but
generically there will less independent invariant functions then
half of the dimension of the phase space. However, the Hamiltonian
flow generated by an adjoint invariant function can be explicitly
described in terms of the factorization on corresponding Lie
group, which is an indication of integrability of such systems. We
will call such systems characteristic ( in algebraic case the
integrals are characters of finite dimensional representations).

As it was already mentioned that the phase space of the Toda system
corresponding to a simple Lie algebra $\mathfrak{g}$ is a special
coadjoint orbit of the Poisson Lie group whose tangent Lie
bialgebra is dual to $\mathfrak{g}$. The integrals are given by
restriction to this orbit of functions which are invariant with
respect to the adjoint action of $\mathfrak{g}$. It is natural to
study Hamiltonian systems generated by such functions on other
orbits. This has been done in \cite{Deift}
\cite{Flash}\cite{Gekht}. In \cite{Deift} the integrability of
such systems on generic orbit has been proven for classical Lie
algebras. In \cite{Flash} it was argued that such systems are
integrable for generic triangular orbits. In \cite{Gekht} it was
proven that such systems on generic orbit are integrable for any
simple finite dimensional Lie algebra $\mathfrak{g}$. An important
result of \cite{Gekht} is that these systems are not only
integrable in the usual Lioville sense but there are also
so-called degenerate integrable systems (with the dimension of the
invariant tori less then half of the dimension of the phase
space). Such degenerate integrable systems were systematically
introduced in \cite{N} and they are also known as systems with
non-commutative integrability (see \cite{Fomenko} where they were
used in a special case). This notion is a classical version of
hidden symmetries in quantum mechanics (see \cite{Pauli} \cite{Winter}).

Similar question can be asked about corresponding Poisson Lie
groups. In this case special symplectic leaves give
generalizations of Toda systems \cite{HKKR} known for $SL_n$ as
relativistic Toda systems.
For generic symplectic leaves of $SL_n$ the
integrability of characteristic systems has been proven in \cite{Li97}.

In this paper we will show the degenerate integrability of a
characteristic Hamiltonian system on any symplectic leaf of any
simple Poisson Lie groups. We will focus on complex algebraic case
and will consider the systems corresponding to real forms only
briefly.

The author would like to thank D. Ben-Zvi, E. Frenkel, M.
Gekhtman, C. Procesi, M. Semenov-Tian-Shansky and M. Yakimov for
interesting discussions.

\section{Characteristic Hamiltonian systems of factorizable Poisson Lie groups}

\subsection{Factorizable Lie bialgebras}
Recall that a Lie bialgebra is a pair $(\g, \delta)$ where $\g$ is
a Lie algebra and the linear map $\delta: \g\to \g\wedge \g$ is a
$\g\wedge \g$-valued 1-cocycle on $\g$ which induces Lie algebra
structure on the dual vector space to $\g$.

A Lie bialgebra $(\g,\delta)$is called factorizable if there
exists $r\in \g\otimes \g$ such that:
\begin{itemize}
\item $r+\sigma(r)$ ($\sigma(x\otimes y)=y\otimes x$) is a
nondegenerate element of $\g\otimes \g$ invariant with respect to
the diagonal adjoint $\g$-action,
\item $\delta(x)=[r,x]$, where the bracket is the diagonal
adjoint action of $x$ on $\g\otimes \g$ and the result is
in $\g\wedge \g \subset \g\o \g$,
\item
\[
[r_{12},r_{13}]+ [r_{12},r_{23}]+[r_{13},r_{23}]=0
\]
where $\g^{\otimes 3} \subset U(\g)^{\otimes 3}$ and
$r_{12}=r\otimes 1, r_{23}=1\otimes r$ etc..
\end{itemize}

Let $\g$ be a factorizable Lie bialgebra with classical r-matrix $r\in
\g\otimes \g$.
\begin{lemma} The linear maps $r_\pm:\g^*\longrightarrow\g$
\[ r_+(l):=\langle r,l\otimes id\rangle,\quad r_-(l):=-\langle r,id\otimes l\rangle\qquad\forall\,l\in\g^*; \]
are Lie bialgebra homomorphisms.
\end{lemma}
The proof of this lemma we will leave as an exercise. It follows
from the classical Yang-Baxter equation for $r$.

The linear map $I:\g^*\to \g, \ l\mapsto r_+(l)-r_-(l)$ is
a linear isomorphism. It is called \emph{the factorization map}.

\begin{corollary} The subspaces $\g_{\pm}={\mbox{Im}}(r_{\pm})$
are Lie subbialgebras in
the Lie bialgebra $\g$.
\end{corollary}

\bigskip

\begin{lemma} {\rm{(1)}} The subspaces
$\N_{\pm}=r_{\pm}({\mbox{ker}}(r_{\mp}))$ are Lie ideals in
$\g_{\pm}$ respectively.

{\rm{(2)}} The map $\theta: \g_+/\N_+\to \g_-/\N_-$ which sends
the residue class of $r_+(e){\mbox{mod}} \ \N_+$ to
$r_-(e){\mbox{mod}} \ \N_-$ is defined and is an isomorphism of
Lie algebras.
\end{lemma}

{\bf Proof.} The first statement follows from the facts that ker$(r_{\pm})$
are Lie ideals in $\g^*$ and that $r_{\pm}$ are Lie algebra
homomorphisms.

Let $\ell\in\g^*$, \ $n\in{\mbox{ker}}(r_+)$, \ $m\in
{\mbox{ker}}(r_-)$. Consider representative $r_+(\ell)+r_+(m)$ of
$r_+(\ell){\mbox{mod}}(\N_+)$ by linearity
$r_+(\ell)+r_+(m)=r_+(\ell +m)= r_+(\ell +n+m)$. Now, $r_-(n+\ell
+m)=r_-(\ell +n)=r_-(\ell) +r_-(n)$ represents the equivalence
class $r_-(\ell){\mbox{mod}}(\N_-)$. Therefore the map $\theta:
\g_+/\N_+\to \g_-/\N_-$ is defined and it is a linear isomorphism.

We leave as an exercise to prove that it is a Lie algebra homomorphism.

\begin{theorem} {\rm{(1)}} Every element $x\in\g$ admits
unique factorization
\[
x=x_+-x_-
\]
where $x_{\pm}\in \g_{\pm}$, and $\theta(x_+ {\mbox{mod}} \ \N_+)=
x_- {\mbox{mod}} \ \N_-$.

{\rm{(2)}} The Lie algebra $\g^*$ is isomorphic to the
following Lie subalgebra of $\g_+\oplus\g_-$:
\[
\{(x_+,x_-)\in \g_+\oplus\g_- \mid \theta(x_+ {\mbox{mod}} \
\N_+)= x_- {\mbox{mod}} \ \N_-\}
\]

{\rm{(3)}} If we model $\g^*$ as in {\em{(2)}}
the factorization map $I:x\mapsto r_+(x)-r_-(x)$ acts as
$(x_+,x_-)\mapsto x_+-x_-$, where on the right side we consider
$\g_{\pm}$ as Lie subalgebras of $\g$.
\end{theorem}

\subsection{Factorizable Poisson Lie groups}

A Poisson Lie group is a Lie group with the such Poison manifold
structure on it that the multiplication map $G\times G\to G$ is a
Poisson map.

There is a bijection between connected simply-connected finite-
dimensional Poisson Lie groups and finite-dimensional Lie
bialgebras. Each finite-dimensional Lie bialgebra can be
"exponentiated" to a connected simply-connected Poisson Lie group
and conversely, a Poisson Lie structure on a finite-dimensional
Poisson Lie group defines a Lie bialgebra structure on the
corresponding Lie algebra. This Lie bialgebra structure is called
tangent Lie bialgebra to a Poisson Lie group. Thus, we have
Poisson Lie counterparts of Lie bialgebras introduced above:
coboundary, quasitriangular, triangular, factorizable.

If $(G,p)$ is a quasitriangular Poisson Lie group, the Poisson tensor
has the following explicit description:
\[
p(x)={\mbox{Ad}}_x(r)-r\in\wedge^2 TG\simeq \wedge^2  g \ .
\]
Here we trivialized the tangent bundle by left translations.
For the Poisson brackets on a quasitriangular  Poisson Lie group we have:
\[
\{f_1,f_2\}=\langle r, d_lf_1\wedge d_lf_2\rangle -
\langle r, d_rf_1\wedge d_rf_2\rangle
\]
where $d_l$ and $d_r$ are, respectively, left and right differentials
on $G$.

For factorizable  Poisson Lie groups  we have
\begin{itemize}
\item maps $r_{\pm}$ lift to Lie group homomorphisms $r_{\pm}:G^*\to G$.

\item $G^{\pm}={\mbox{Im}}(r_{\pm})\subset G$ are Poisson Lie subgroups
(connected simply connected)
\item  $N^{\pm}={\mbox{Im}}({\mbox{ker}}(r_{\pm}))\subset G^{\pm}$ are
normal Lie subgroups
\item Lie algebra isomorphism $\theta: \g_+/\N^+\stackrel{\sim}{\rightarrow} \g_-/\N_-$ lifts to Lie group isomorphism
$\theta: G^+/N^+ \stackrel{\sim}{\rightarrow} G^-/N^-$.
\item Lie group $G^*$ can be modeled as:
\[
G^*=\{(g_+,g_-)\in G^+\times G^- \mid
\theta(g_+{\mbox{mod}} \ N^+)=g_-{\mbox{mod}} \ N^-\}
\]
\item There exist open dense subsets
$G\,',G\,''\subset G$ such that for each $g\in G\,'$ there exists unique
factorization $g=g_+g_-^{-1}$, $\theta(g_+{\mbox{mod}} \
N_+)=g_-{\mbox{mod}} \ N^-$ and for each $g\in G\,''$ there exists unique
factorization $g=g^{-1}_-g_+$ with the same conditions on $g_{\pm}$.
\item Left, respectively right, {\em factorization maps} \
$G^*\to G$ map $(g_+,g_-)$  to $(g_+g_-^{-1})$, respectively to
$g^{-1}_-g_+$.
\end{itemize}
Here we assume that $G^*$ is represented as a subgroup of
$G^+\times G^-$.

\subsection{The double}The double $D(\g)$ of the Lie bialgebra $\g$
is the Lie bialgebra which  is the direct sum $\g\oplus
{\g^*}^{op}$ of Lie coalgebras. The  Lie bracket on it is
determined uniquely by the requirement that the natural bilinear
form $<(x,l),(y,m)>=l(y)+m(x)$ is $D(\g)$-invariant and the
isotropic subspaces $\g$ and $\g^*$ are Lie subalgebras. We will
denote these Lie bialgebra inclusions $i:\g\to D(\g)$ and $j:
{\g^*}^{op}\to D(\g)$.

The double $D(G)$ of $G$ is the connected, simply connected
Poisson Lie group having $D(\g)$ as its Lie bialgebra. The maps
$i$ and $j$ lift to injective Poisson Lie maps $i:G \to D(G)$,
$j:{G^*}^{op} \to D(G)$ and consequently to a map $\mu\circ
(i\times j):G\times {G^*}^{op} \to D(G)$: $(x,y)\mapsto i(x)j(y)$
which is also a local Poisson isomorphism. Here ${G^*}^{op}$ is
the Lie Poisson group $G^*$ with the opposite Poisson structure.

\subsection{Symplectic leaves in Poisson Lie groups}
The Poisson Lie group ${G^*}^{op}$ acts on $D(G)$ via left
multiplication, $y\cdot x:= j(y)x$. We also have a map $\varphi
:G\to D(G)/j(G^{*op})$ which is the composition of $i$ with the
natural projection. This map is Poisson. In a neighborhood of the
identity this map $\varphi$ is a Poisson isomorphism.  The map
$\varphi$ has open dense range but it is not surgective if the
factorization problem in $D(G)$ does not have global solution. It
intertwines local dressing action of $G^*$ on $G$ \cite{STS} and
the action of $G^*$ on the cosets.

The symplectic leaves of $G$ are orbits of dressing action of
$G^{*op}$. Or, equivalently, they are connected components of preimages of left
$G^{*op}$-orbits in $D(G)/j(G^{*op})$. Notice that this
description does use the Poisson structure on the groups and therefore
we can use notation $G^*$ instead of $G^{*op}$ without a
danger of confusion.

\subsection{Characteristic Hamiltonian systems on factorizable Poisson Lie groups}

Let $(G,p)$ be a factorizable Poisson-Lie group. Let $I(G)\subset C^\infty(G)$ be the subspace of $Ad_G$-invariant functions on $G$.

\begin{theorem} \label{main}\hfill
\begin{enumerate}
\item[i)] $I(G)$ is a commutative Poisson algebra in $C^\infty(G)$.
\item[ii)] In a neighborhood of $t=0$ the flow lines of the Hamiltonian flow induced by $H\in I(G)$ passing through $x\in G$ at $t=0$ have the form
\[ x(t)=g_\pm(t)^{-1}x g_\pm(t), \]
where the mappings $g_\pm(t)$ are determined by
\[ g_+(t)g_-(t)^{-1}=\exp\left(tI\left(d H(x)\right)\right), \]
and $I:\g^*\longrightarrow\g$ is the factorization isomorphism.
\end{enumerate}
\end{theorem}

\begin{defin} Characteristic Hamiltonian system on a factorizable
Poisson Lie group is a Hamiltonian system whose phase space
is a symplectic leaf of a factorizable Poisson Lie group
and whose Hamiltonian is an adjoint invariant function on $G$.
\end{defin}

The theorem above implies that the equations of motion of a
characteristic Hamiltonian system on a factorizable Poisson Lie
group can be solved via factorization.  Below we will show that
such systems are integrable for all simple Poisson Lie groups.

\section{Standard Poisson Lie structure on simple Lie groups}

\subsection{Standard factorizable Lie bi-algebra structure on simple Lie algebras
and standard Poisson Lie structure on simple Lie groups}
Let $\g$ be a simple complex Lie algebra. Fix a Borel subalgebra $\borel$.
Let $(H_i,e_i,f_i)$ $i=1,\ldots,r=\mbox{rank}(\g)$
be elements of the Chevalley basis of $\g$ for this choice of Borel subalgebra
which correspond to simple roots. It is well known that $\g$
is freely generated by $(H_i,e_i,f_i)$ modulo determining relations
\begin{eqnarray}
[H_i,H_j] &=& 0, \label{rel-1}\\[0pt]
[H_i,e_j] &=& a_{ij}e_j, \\[0pt]
[H_i,f_j] &=& -a_{ij}f_j, \\[0pt]
[e_i,f_j] &=& \delta_{ij}H_i, \\
(ad_{e_i})^{1-a_{ij}}(e_j) &=& 0,\qquad i\neq j, \\
(ad_{f_i})^{1-a_{ij}}(f_j) &=& 0,\qquad i\neq j, \label{rel-6}
\end{eqnarray}
where $(a_{ij})$ denotes the Cartan matrix of $\g$.\\

\noindent

Consider a linear map $\delta:\g\rightarrow\g\wedge\g$  acting on the
generators as
\begin{eqnarray}
\delta(H_i) &=& 0, \\
\delta(e_i) &=& \frac 12 d_i H_i\wedge e_i, \\
\delta(f_i) &=& \frac 12 d_i H_i\wedge f_i,
\end{eqnarray}
where $d_i$ is the length of the i-th root, in particular,
\[
d_i a_{ij}=a_{ji}d_j.
\]

\begin{theorem}
There exists unique such linear map $\delta:\g\rightarrow\g\wedge\g$
which has the following properties:
\begin{enumerate}
\item[i)] $\delta$ is a 1-cocycle.
\item[ii)] $(\delta\wedge id)\circ\delta=0$.
\end{enumerate}
\end{theorem}

\begin{rem} \rm
The definition of the standard Lie bialgebra structure on $\g$ requires the choice of a Borel subalgebra $\borel\subset\g$.
\end{rem}

\begin{theorem}
$(\g,\delta)$ is factorizable with
\[ r=\frac{1}{2}\sum_{ij}(B^{-1})_{ij}H_i\otimes H_j+\sum_{a>0}e_\alpha\otimes f_\alpha,\]
where $B_{ij}=d_i a_{ij}$ is the symmetrized Cartan matrix.
\end{theorem}

This induces the Poisson Lie structure on $G$ for which the Lie
bialgebra described above is the tangent Lie bialgebra. The Borel
subgroup $B$ and its opposite $B^-$ are Poisson Lie subgroups.

The Lie bialgebra $Lie(G)$ is isomorphic to the double of the Lie
bialgebra $Lie(B)$ quotioned by the diagonally embedded Cartan
subalgebra \cite{Drinfeld87}.

We denote by $\N^+$ and $N^-$ the nilpotent subgroups of $B^-$ and
$B^-$, respectively. Since $H=B^+\cap B^-$ we have two natural
projections $[ \ ]_0:B^+\to B^+/N^+\cong H$ $b\mapsto [b]_0$ and
$[ \ ]_0:B^-\to B^-/N^-\cong H$. We shall also write $B^+$ and
$N^+$ for $B$ and $N$, respectively.

\subsection{Symplectic leaves}
\subsubsection{Symplectic leaves of $B^{\pm}$.} It is known that
$(B^+)^{*{\mb{op}}}\simeq B^-$ as a Poisson Lie group and that
$D(B^+)\simeq G\x H$ as a Lie group. Fix these isomorphisms. The
double $D(B^+)$ is a factorizable Poisson Lie group with Poisson
Lie imbeddings $i:B^+\hookrightarrow D(B^+)$, \
$j:B^-\hookrightarrow D(B^+)$
\[
i(b^+)=(b^+,[b^+]_0) \ , \qquad j(b^-)=(b^-,[b^-]_0^{-1})
\]
The group $B^-$ acts on cosets $D(B^+)/j(B^-)$ by multiplication
from the left. Define Lie subalgebras ${\mathfrak h}^w=coker\{w-id\}$ and
${\mathfrak h}_w=ker\{w-id\}$ of the Cartan subalgebra. Here the element $w\in W$
of the Weyl group $W$ of $G$ is considered as a linear operator
on the Cartan subalgebra. Let $H^w$ and $H_w$ be corresponding Lie subgroups
in $H$.

We have left Bruhat decoposition of $D(B^+)$:
\[
D(B^+)=\sqcup_{w\in W} D(B^+)_w, \ \ D(B^+)_w=B^-wB^-\times H.
\]

Orbits of the action of $B^-$ on the cosets $D(B^+)/j(B^-)$ have
the following structure:
\begin{itemize}
\item $j(B^-)\backslash D(B^+)_w/j(B^-)\cong H_w$
\item each $B^-$ orbit in $D(B^+)_w/j(B^-)$ is isomorphic to
$N^-_w\times H^w$.
\end{itemize}
Here $N^-_w$ is the subspace of the nilpotent subgroup
$N^-$ generated by one parametric subgroups generated
by those negative roots which remain negative after the
action of $w$.

Symplectic leaves of $B^+$ are irreducible components of preimages
of $B^-$-orbits with respect to the map
\[
\varphi: B^+\stackrel{i}{\hookrightarrow} D(B^+)\longrightarrow
D(B^+)/j(B^-) \ .
\]

The image of $\varphi$ intersects $j(B^-)$ orbits. Consider sets
$B^+_w=B^+\cap B^-wB^-$. It is a Poisson subvariety \cite{DCP}
\cite{HKKR} in $G$.

\begin{lemma}
Let ${\mathcal O}_w\subset B^-wB^-\times H/j(B^-)$ be an orbit of
the right $B^-$ action on this coset, then $\varphi(B^+_w)\cap
{\mathcal O}_w\subset{\mathcal O}_w$ is Zariski open.
\end{lemma}

From this and from the fact that $\varphi$ is a cover map
one can show that $B^+_w=B^+\cap B^-wB^-$ is a Poisson
subvariety in $B^+$ with symplectic leaves of have dimension
$\ell(w)+$ corank $(w-{\mb{id}})$ \cite{DCP}.

One can give "explicit" description of symplectic
leaves as connected components of Casimir functions.

According to \cite{FZ} define generalized minors as the following
functions on the group $G$. Let $G_0$ be the subset in $G$ formed
by elements who have Gaussian factorization $x=[x]_-[x]_0[x]_+$
with $[x]_{\pm} \in N^{\pm}$ and $[x]_0\in H$. For a weight
$\lambda$ define function $$\D_\lambda(x)=[x]_0^\lambda$$.

Let $\omega_i$ be a highest weight of $i$-th fundamental representation
of $G$ and $\bar{\bar{u}}$ and $\bar{v}$ are special representatives
of elements $u,v \in W$ in $G$.  Generalized minors are the following
functions:
\begin{equation}\label{DDelta}
\D_{u\omega_i,v\omega_i}=\D_{\omega_i}(\bar{\bar{u}}^{-1}x\bar{v})
\end{equation}
where $\bar{u}$ and $\bar{\bar{u}}$ are special representatives of
the element $u\in W$ (see \cite{FZ} for details).

\begin{lemma} The generalized minors $\D_{\O_i,w^{-1}\O_i}(x)$
do not vanish on $B^+_w$ and generalized minors
$\D_{w\O_i,\O_i}(x)$ do not vanish on $B^-_w$.
\end{lemma}
{\it Proof}. Let us prove that $\D_{\omega_i,w^{-1}\omega_i}$ does
not vanish on $B^-wB^-$. Consider $x=b_-\overline{w}\beta_-\in
B^-wB^-$ and let $b_-=[b_-]_0[b_-]_-$ and
$\beta_-=[\beta_-]_0[\beta_-]_-$ be Gauss decompositions, then we
have
\begin{eqnarray*}
\D_{\omega_i,w^{-1}\omega_i}(b_-\overline{w}\beta_-)&=&
\D_{\omega_i}([b_-]_0\overline{w}\beta_-\overline{w}^{-1})=\\
\D_{\omega_i,\omega_i}([b_-]_0w([\beta_-]_0)
\overline{w}[\beta_-]_-\overline{w}^{-1})&=&
\D_{\omega_i}([b_-]_0w([\beta_-]_0))
\end{eqnarray*}
The last function does not vanish, which proves the first statement
of the lemma. The proof of the last one is similar.

Here and below we will denote by $\Lambda=\z\omega_1\oplus\dots\z\omega_r$
the weight lattice in ${\mathfrak h}^*$.
\begin{proposition} Symplectic leaves of $B^+_w$ are irreducible
components of level surfaces of functions
\begin{equation} \label{casimirs-b}
c^+_{w,t}(x)=\prod^r_{i=1}\D_{\O_i,w^{-1}\O_i}(x)^{t_i}
([x]^{\O_i}_0)^{t_i}
\end{equation}
where $t=\sum^r_{i=1}t_i\O_i\in{\rm{ker}}_{\Lambda}(w-{\rm{id}})\in
{\mathfrak h}^*$.
\end{proposition}

{\bf Proof.} First, let us prove that functions $c^+_{w,t}(x)$ are
invariant with respect to the (local) dressing action.

Dressing action of $b_-\in B^-$ on $b^+\in B^+$ is the map
$b_-:b_+\mapsto b^{b_-}_+$ given by the solution to the
factorization equations:
\begin{eqnarray*}
b_+b_-^{-1}&=& (b_-^{b_+})^{-1} b_+^{b_-} \\ [0.125in]
\left[b_+\right]_0  \ [b_-]_0 &=& [(b_-^{b_+}]_0 \ [(b_-^{b_-}]_0
\end{eqnarray*}
where  $b_{\mp}^{b_{\pm}}\in B^{\mp}$. This system has a unique
solution when $b_-$ is sufficiently close to 1.

1. Assume that $[b_-]_-=1$, then $b_-=[b_-]_0\in H$ and
\[
b_+^{b_-}=[b_-]_0 \ b_+[b_-]^{-1}_0 \ .
\]
It is clear that functions $c^+_{u,t}$ are invariant with respect
to such action of $H\subset B^-$ \ iff \ $t=u(t)$.

\bs 2. Assume that $[b_-]_0=1$, then  $b_-\in N^-$. Denote by
${\tilde x}_{\pm}\in B^{\pm}$ the result of ``opposite"
factorization of $(x,h)\in G\x H=D(B^+)$:
\[
{\tilde x}_-^{-1}{\tilde x}_+=x \ , \qquad [{\tilde x}_-]_0 \
[{\tilde x}_+]_0 =h \ .
\]
Then we have
\begin{eqnarray}\label{factorb}
b_+b_-^{-1} &=& b_-^{-1}(b_-b_+b_-^{-1}) = b_-^{-1}
(\widetilde{b_-b_+b_-^{-1}})^{-1}_- \
(\widetilde{b_-b_+b_-^{-1}})_+ \nonumber \\ [0.25in]
\left[b_+\right]_0 &=& [(\widetilde{b_-b_+b_-^{-1}})_-]_0 \
[(\widetilde{b_-b_+b_-^{-1}})_+]_0
\end{eqnarray}

\begin{lemma} \label{Delta}
Functions $\D_{\O_i,u^{-1}\O_i}(x)$ on $B_-uB_-\subset G$ are
invariant with respect to the left action of $N^-$.
\end{lemma}

{\bf Proof.} Let $x\in B_-uB_-$ and  $x=\beta'_-{\bar u}
\beta''_-$ where ${\bar u}\in G$ is a special representative of
$u\in W$ (see {FZ} for definition of ${\bar u}$), $\beta'_-\in
N^-$ and $\beta''_-\in B^-$. Then for $n_-\in N^-$ we have
\begin{eqnarray*}
 \D_{\O_i,u^{-1}\O_i}(xn_-) &=&
\D_{\O_i}(xn_- {\bar u}^{-1}) = \D_{\O_i}(\beta'_-{\bar u}
\beta''_- n_-{\bar u}^{-1})  \nonumber\\ [0.125in]
 &=& \D_{\O_i}(u([\beta''_-]_0){\bar u}\beta''_-n_-{\bar
u}^{-1})
\end{eqnarray*}
The element ${\bar u}\beta''_-n_-{\bar u}^{-1}$ always admits
factorization into the product $x_-x_+$, \ $x_{\mp}\in N^{\mp}$
and therefore
\[
\D_{\O_i,u^{-1}\O_i}(xn_-)= \D_{\O_i}(u[\beta'']_0) =
\D_{\O_i,u^{-1}\O_i}(x) \ .
\]
Now let us compute \ $\D_{\O_i,u^{-1}\O_i}(b_+^{b_-})$ for $b_-\in
N^-$:
\begin{eqnarray*}
\D_{\O_i,u^{-1}\O_i}(b_+^{b_-})&=&
\D_{\O_i,u^{-1}\O_i}(\widetilde{b_-b_+b_-^{-1}})_+ \nonumber \\
[0.125in] &=& \D_{\O_i,u^{-1}\O_i}(b_-b_+b_-^{-1})
([b_-b_+b_-^{-1})_-]^{\O_i}_0)^{-1} \nonumber \\ [0.125in] &=&
[b_+]_0^{\O_i} ([(\widetilde{b_-b_+b_-^{-1}})_+]_0^{\O_i})^{-1}
\D_{\O_i,u^{-1}\O_i}(b_+b_-^{-1})\\ [0.125in] &=& [b_+]_0^{\O_i}
([b_+^{b_-}]_0^{\O_i})^{-1} \D_{\O_i,u^{-1}\O_i}(b_+) \ .
\nonumber
\end{eqnarray*}
Here we used (\ref{factorb}) (factor b) and the Lemma. Thus, the
function
\[
\D_{\O_i,u^{-1}\O_i}(b_+) \ [b_+]_0^{\O_i}
\]
is invariant with respect to the dressing action of $N^-$.

Thus, functions (\ref{casimirs-b} ) are invariant with respect to
the action of $H$ and $N^-$ and therefore they are invariant with
respect to the dressing action of $B^-$ on $B^-uB^-$ and therefore
they are Poisson Casimirs for $B^-uB^-$. They do not vanish on $B^-uB^-$
and therefore
their level surfaces are Poisson subvarieties and they form a fiber
bundle over $({\mathbb C}^{\times})^{\mbox{\footnotesize
corank}(u-\mb{\footnotesize id})}$. Dimension of fibers coincide
with the dimension of symplectic leaves of $B^+$ which are in
$B^+\cap B^-uB^-$ which proves the proposition.

Similarly for $B^-$. The subsets $B^-_w=B^-\cap B^+wb^+$ are
Poisson subvarieties whose symplectic leaves have
dimension $\ell(w)+$ corank$(w-{\mb{id}})$.

\begin{proposition} Symplectic leaves of $B^-_w$ are irreducible
components of level surfaces of
functions
\[
c^-_{w,s}(x)=\prod^r_{i=1}\D_{w\O_i,\O_i}(x)^{-s_i}
([x]^{\O_i}_0)^{s_i}
\]
where $s=\sum^r_{i=1}s_i\O_i\in{\rm{ker}}_{\Lambda}(w-{\rm{id}})\subset
{\mathfrak h}^*$.
\end{proposition}

\subsubsection{Symplectic leaves of $D(B^+)$.} The dual Poisson Lie
group to $D(B^+)=G\x H$ can be identified with $B^+\x B^-$ (as a
Lie group). We also have an isomorphism of Lie groups
$D(D(B^+))\simeq D(B^+)\x D(B^+)$. Fix these isomorphisms.

The Poisson Lie imbeddings $i:D(B^+)\hookrightarrow D(B^+)\x
D(B^+)$, $j:D(B_+)^{*{\mb{op}}} \hookrightarrow D(B^+)\x D(B^+)$
are
\[
i(g,h)=((g,h),(g,h)),\ \ j((b_+,b_-))=
((b_+,[b_+]_0),(b_-,[b_-]^{-1}_0) \ .
\]
Symplectic leaves of $D(B^+)=G\x H$ are connected components of
preimages of left $j(D(B^+)^{*{\mb{op}}}$-orbits in $D(B^+)\x
D(B^+)/jD(B^+)^{*{\mb{op}}}$ with respect to the map
\[
\phi:D(B^+)\hookrightarrow D(B^+)^{\times 2}\longrightarrow
D(B^+)^{\times 2}/j(D(B^+)^{*{\mb{op}}}) \ .
\]
Let $G^{u,v}=B^+uB^+\cap B^-vB^-$ be the double Bruhat cell
corresponding to the pair $(u,v)\in W\x W$. One can show that
$G^{u,v}\x H$ is a Poisson subvariety with symplectic leaves of
dimension
$\ell(u)+\ell(v)+{\mb{corank}}(u-{\mb{id}})+{\mb{corank}}(v-{\mb{id}})$
\cite{HKKR}. Since generalized minors $\D_{\O_i,u^{-1}\O_i}(x)$ do
not vanish on $B^+_u$ and generalized minors $\D_{v\O_i,\O_i}(x)$
do not vanish on $B^-_v$ neither of them vanish on the
intersection $G^{u,v}=B^+uB^+\cap B^-vB^-$.

Next proposition describes symplectic leaves of $D(B^+)$ in terms of level
sets of Casimir functions.

\begin{proposition}  Symplectic leaves of $G^{u,v}\x H$ are irreducible
components of level surfaces of functions
\[
c_{u,v,s,t}(x,h)=c^+_{u,s}(x,h)c^-_{v,t}(x,h)
\]
where \ $s=\sum^t_{i=1}s_i\O_i\in{\mb{ker}}_{\Lambda}(u-{\mb{id}})\subset
{\mathfrak h}^*$, \
$t=\sum^r_{i=1}t_i\O_i\in{\mb{ker}}_{\Lambda}(v-{\mb{id}})\subset {\mathfrak
h} ^*$ and
\begin{eqnarray*}
c^+_{u,s}(x,h) &=& \prod^r_{i=1} \D_{\O_i,u^{-1}\O_i}(x)^{s_i}
(h^{\O_i})^{s_i} \ , \\ c^-_{v,t}(x,h) &=& \prod^r_{i=1}
\D_{v\O_i,\O_i}(x)^{t_i} (h^{\O_i})^{-t_i} \ .
\end{eqnarray*}
\end{proposition}

{\bf Proof.} We have the following natural identifications of
groups:
\[
D(B^+) = G\x H \ , \qquad D(B^+)^*=B^+\x B^- \ , \qquad
D(D(B^+))=D(B^+)\x D(B^+) \ .
\]
Left and right factorizations of an element in $D(D(B^+))$ are:
\[
((g_1,h_1),(g_2,h_2)) = ((g,h)(g,h))(i(\xi_+),j(\xi_-)) =
(i(\xii_+),j(\xii_-))((\Ge,\h),(\Ge,\h)) \ ,
\]
or, in components:
\[
g\xi_{\pm} =\xii_{\pm}\Ge \ , \qquad h[\xi_+]_0=[\xii_{\pm}]^{\pm
1}_0\h \ .
\]
Such factorization exists on an open dense subset of $D(D(B))$.

The map $\xi:D(B^+)\to D(B^+)$ acting as $\xi:g\mapsto\Ge$, \
$h\to\h$ determines the (local) dressing action. For each $(g,h)$,
it is defined for $\xi$ sufficiently close to 1.

Using Lemma \ref{Delta} we obtain
\[
\D_{\O_i,u^{-1}\O_i}(\Ge) =
\D_{\O_i,u^{-1}\O_i}(\xii_-^{-1}g\xi_-) =[\xii_-]_0^{-\O_i}
[\xi_-]_0^{u^{-1}(\O_i)} \D_{\O_i,u^{-1}\O_i}(g) \ .
\]
Similarly,
\[
\D_{v\O_{i},\O_i}(\Ge) = \D_{v\O_i,\O_i}(\xii_-^{-1}g\xi_+)
=[\xi_+]_0^{\O_i} [\xii_+]_0^{-v^{-1}(\O_i)} \D_{v\O_i,\O_i}(g) \
.
\]
From here we see that
\begin{eqnarray*}
&&\prod^r_{i,j=1} \D_{\O_i,u^{-1}\O_i}(\Ge)^{t_i}
\D_{v\O_j,\O_j}(\Ge)^{s_j} \\ &&\quad
=\prod^r_{i,j=1}\D_{\O_i,u^{-1}\O_i}(g)^{t_i}
\D_{v\O_j,\O_j}(g)^{s_j}\\ &&\quad \dt [\xi_-]_0^{-t+u^{-1}(t)}
[\xi_+]_0^{-v^{-1}(s)+s} (h\h^{-1})^{-v^{-1}(s)+t}
\end{eqnarray*}
Thus, the functions $c_{u,v,t,s}$ on $D(B^+)$ are invariant with
respect to dressing transformations if and only if
\[
t=u(t) \ , \qquad s=v(s) \ .
\]
Thus, functions $c_{u,v,t,s}$ are Casimirs in the Poisson algebras
of functions on $G^{u,v}\times H$. On the other hand, they do not vanish
on $G^{u,v}\times H$. Therefore their level surfaces form a fiber bundle
over the torus $({\mathbb C}^{\times})^{\mbox{\footnotesize
corank}(uv^{-1}-\mb{\footnotesize id})}\x ({\mathbb
C}^{\times})^{\mbox{\footnotesize corank}(u-\mb{\footnotesize
id})}$.
The fibers are Poisson subvarieties and they  have the same dimension as symplectic leaves of
$G^{u,v}\times H$.

\subsubsection{ Symplectic leaves of $G$.}

Since $G$ is a factorizable Poisson Lie group and its symplectic
leaves can be described very similarly to those of $D(B^+)$
\cite{HL}. Taking into account that $D(G)=G\times G$ we have the
composition of Poisson maps:
\[
\phi: G\hookrightarrow G\times G \to (G\times G)/j(G^*)
\]
where $j(G^*)=\{(b,b_-)\in B^+\times B^-| [b]=[b_-]^{-1}\}$.
Connected components of primages of $j(G^*)$-orbits on cosets
are symplectic leaves of $G$.

The following proposition describes symplectic leaves of $G$ in
terms of Casimir functions.

\begin{proposition}
Double Bruhat cell $G^{u,v}=B^+uB^+\cap B^-vB^-$ is a Poisson
subvariety in $G$ with symplectic leaves of dimension
$\ell(u)+\ell(v)+{\rm{corank}}(uv^{-1}-{\rm{id}})$. They
are irreducible components of level surfaces of functions
\[
c_{u,v,t}(x)=\prod^r_{i=1}\D_{v\O_i,\O_i}(x)^{t_i}
\D_{\O_i,u^{-1}\O_i}(x)^{u^{-1}(t)_i}
\]
where \ $t=\sum^r_{i=1}t_i\O_i\in{\rm{ker}}_{\Lambda}(uv^{-1}-{\rm{id}})\subset
{\mathfrak h}^*$.
\end{proposition}
{\bf Proof.} The (local) dressing action of
$G^*=\{(\xi_+,\xi_-)\mid [\xi_+]_0=[\xi_-]^{-1}_0\}\subset B^+\x
B^-$ on $G$ is given by the map $(\xi_+,\xi_-):g\mapsto\Ge$, where
$\Ge$ is the solution to the factorization problem:
\begin{eqnarray}\label{xi}
g\xi_{\pm} &=& \xii_{\pm}\Ge \nonumber\\ \left[\xi_+\right]_0 &=&
[\xi_-]^{-1}_0 \ , \qquad [\xii_+]_0=[\xii_-]_0^{-1} \ .
\end{eqnarray}
Such $\Ge,\xii_{\pm}$ exist for each $g$ where $(\xi_+,\xi_-)$ are
sufficiently close to 1. Double Bruhat cells are invariant
submanifolds for this action.

On $G^{u,v}$ we have (as above for $D(B^+)$)
\begin{eqnarray*}
\D_{\O_i,u^{-1}\O_i}(\Ge) &=& [\xii_-]_0^{-\O_i}
[\xi_-]_0^{u^{-1}(\O_i)} \D_{\O_i,u^{-1}\O_i}(g) \\ \D_
{v\O_{\ell},\O_{\ell}}(\Ge) &=& [\xi_+]_0^{\O_i}
[\xii_+]_0^{-v^{-1}(\O_i)} \D_{v\O_i,\O_i}(g)
\end{eqnarray*}
Therefore
\begin{eqnarray}\label{Del-Cas}
&&\prod^r_{ij=1} \D_{\O_i,u^{-1}\O_i}(\Ge)^{t_i}
\D_{v\O_j,\O_j}(\Ge)^{s_j} \\ && \quad = \prod^r_{ij=1}
\D_{\O_i,u^{-1}\O_i}(g)^{t_i} \D_{v\O_j,\O_j}(g)^{s_j} \\ &&\quad
\dt [\xii_-]^{-t}_0 [\xi_-]_0^{u^{-1}(t)} [\xii_+]_0^{-v^{-1}(s)}
[\xi_+]_0^s
\end{eqnarray}
Using relations (\ref{xi}) we can write  the ``$\xi$"-factor as
\[
[\xii_-]_0^{-t+v(s)} [\xi_-]_0^{u^{-1}(t)-s} \ .
\]
Thus, functions (\ref{Del-Cas}) are invariant if and only if
\[
t=(s) \ , \qquad s=u^{-1}(t) \ .
\]
Thus, the functions (\ref{Del-Cas}) are Poisson Casimirs and they
do not vanish on $G^{u,v}$. Therefore their level sets form the fiber bundle over
\[
({\mathbb C}^{\times})^{\mbox{corank}(uv^{-1}-\mb{id})}
\]
with fibers being Poisson submanifolds. The dimension of the
fibers is the same as of the symplectic leaves of $G^{u,v}$. This
proves the theorem.

\bs

\section{Degenerate integrability of Hamiltonian systems}
\label{non-com}

 The notion of degenerate integrability was
introduced in \cite{N}. First examples of such systems were
known long before ( see for example \cite{Pauli} \cite{Winter}) with the
model of the hydrogen atom \cite{Pauli} as a classical example.
Such systems are also known as superintegrable systems \cite{Winter}
and as systems with non-commutative integrability \cite{Fomenko}.

We will say that a subalgebra $A$ of the algebra of functions on a
smooth manifold  $\M$ is generated by functions $f_1,\dots,f_n$ if
for each function $f\in A$ the form $df\wedge\wedge
f_1\wedge\dots\wedge\f_n$.

Assume that we have the following structure on a real symplectic
manifold $(\M_{2n}, \omega)$.

\begin{itemize}
\item $2n-k$ independent functions $J_1, \dots, J_{2n-k}$
generating Poisson subalgebra $C_J(\M)$ in $C(\M)$.

\item $k$ independent functions $I_1, \dots, I_k$ generating
Poisson center of the Poisson subalgebra $C_J(\M)$.
\end{itemize}

Let $H\in C(\M)$ be a function which Poisson commute with
$J_1,\dots, J_{2n-k}$:
\[
\{H,J_i\}=0, \ \ i=1,\dots, 2n-k
\]

We will say that the level surface $\M(c_1,\dots,c_{2n-k})=\{x\in
\M | J_i(x)=c_i\}$ of functions $J_i$ is called generic relative
to functions $I_1,\dots,I_n$ if the form $dI_1\wedge\dots\wedge
dI_{k}$ does not vanish identically on it. Then the following is
true \cite{N}:

\begin{theorem}\label{deg}
\begin{enumerate}
\item Flow lines of $H$ are parallel to level surfaces of $J_i$.
\item Each connected component of a generic level surface
has canonical affine structure generated by the flow lines of
$I_1,\dots, I_k$.
\item The flow lines of $H$ are linear in this affine
structure.

\end{enumerate}
\end{theorem}

When $k=n$ this theorem reduces to the Liouville integrability
\cite{Arnold89}.

When $\M_2n=N_{2k}\times \tilde{N}_{2n-2k}$ and functions
$I_1,\dots, I_k$ are constant along $\tilde{N}_{2n-2k}$ the
degenerate integrability of a system with commuting integrals
$I_1,\dots, I_k$ is equivalent to the Liouville integrability of
the system on $N_k$ with integrals $I_i|_N$ .

Because the theorem \ref{deg} generalizes the Liouville theorem to
the case when the dimension of invariant tori is less the the
dimension of $\M$ we will call these system degenerate integrable
systems.

Geometrically, the structure described above means that we have
two Poisson projections
\begin{equation}\label{psi-pi}
\M_{2n}\stk{\psi}{\longrightarrow} B^J_{2n-k}
\stk{\pi}\longrightarrow B^I_{k}
\end{equation}
where $B^J$ and $B^I$ are Poisson manifolds (level surfaces of
$J_i$ and $I_i$ respectively) and $B^I$ has  the trivial Poisson
structure .

Degenerate integrable systems admit action-angle variables. Let us
call the point $b\in B_I$ {\it regular} if the connected
components of fibers of the preimage ${\pi}^{-1}(b)$ are Poisson
submanifolds in $B_J$ which consist of a single open dense
symplectic leaf.

Let $a\in B^I$ be a regular point and $D$ be  an open neighborhood
of $a$. Choose a generic  point $c\in B_J$ which belongs to an
open dense symplectic leaf of one of the connected components of
${\pi}^{-1}(b)$. Let $U$ be a neighborhood of $c$. Choose the
trivialization of ${\pi}$ over $D$:
\[
{\pi}^{-1}(D)\simeq {\pi}^{-1}(b)\x D \ .
\]
Let $\tilde{U}\subset{\pi}^{-1}(D)$ be a neighborhood of $c$
such that with respect to this trivialization,
\[
\tilde{U}\simeq U\x D \ .
\]
Together with the choice of the trivialization of $\pi$ over $W$
this gives the isomorphism
\[
f:{\psi}^{-1}(\tilde{U})\simeq{\psi}^{-1}(c)\x U\x D \ .
\]

The functions $I_1,\dots, I_k$ give a local coordinate system on
$D$. Their Hamiltonian flows generate $k$ independent Hamiltonian
vector fields on ${\psi}^{-1}$. Define affine coordinates
$\phi_1, \dots, \phi_k$ on this level surface as natural
coordinates along these vector fields.

\begin{theorem} Assume that ${\psi}^{-1}(c)$ is compact. Then there exist
a trivialization $f:{\psi}^{-1}(\tilde{U})\simeq{\psi}^{-1}(c)\x
U\x D$ such that the symplectic form $\omega$ on ${\mathcal M}$
has the form
\[
\omega= f^*\left( \sum^k_{i=1} dI_i\wedge d\varphi_i +
\pi^*(\omega_{{\pi}^{-1}(b)})\right)
\]
where $\omega_{{\pi}^{-1}(b)}$ is the symplectic form on the
open dense symplectic leaf of the connected component of
${\pi}^{-1}(b)$ which contains $c$.
\end{theorem}

 The coordinates $\phi_i, I_i$ are called action-angle variable
 for degenerate integrable systems.

 One can replace real smooth manifolds by complex
manifolds (complex algebraic) and Poisson structures by complex
holomorphic (complex algebraic) structures. In this paper we
assume that $\mathcal M$ is an affine algebraic variety. Poisson
structure on $\mathcal M$ is determined by the structure of
Poisson algebra on the ring of functions.

A degenerate integrable system on an algebraic symplectic manifold
$\mathcal M$ consists of Poisson subalgebra $J$ in the algebra of
functions on $\mathcal M$ with $dim(Spec(J))=2n-k$ and if $Z(J)$
is the Poisson center of $J$, then $dim(Spec(Z(J)))=k$.
Here
$Spec(A)$ is the spectrum of primitive ideals of $A$.

\section{Degenerate integrability of characteristic systems on standard
simple Lie groups with Poisson Lie structure }

\subsection{Poisson structure on $G\times G//{Ad_{G^*}}$}

\subsubsection{The map $\psi$}

Let $M$ be an affine algebraic variety and let $G$ be an algebraic
Lie groups acting on it. Denote by $C_G(M)$ the algebra of
$G$-invariant functions on $M$. We will use notation $M//G$ for
the affine variety which is the categorical quotient,
$M//G=Spec(C_G(M))$. Assume that $M$ is Poisson (the ring of
functions on $M$ is a Poisson algebra). If $M$ has a Poisson
structure, $G$ is a Poisson Lie group and the action is Poisson,
then the algebra $C_G(M)$ is a Poisson subalgebra in $C(G)$ and
therefore $M//G=Spec(C_G(M))$ has natural Poisson structure on it.

\begin{theorem}\label{quot} Let $G$ be an algberaic Poisson Lie group.
The projection $D(G)\to D(G)//Ad_{G^*}$ is a Poisson map.
\end{theorem}

Proof. Similar to the previous section consider two functions $f$
and $g$ on $D(G)$ which are invariant with respect to the adjoint
action of the subgroup $G^*$. Let $x\in D(G)$ and $b\in G^*$, then
\begin{eqnarray}
\{f,g\}(bxb^{-1})&=&<r_D,d_+f\otimes d_-g-d'_+f\otimes
d'_-g>(bxb^{-1})\\
&=&<(Ad_b\otimes Ad_b)(r_D),(df-d'f)\otimes d_-g>(x)\\
&=&<(Ad^*_b\otimes Ad_b)(r_D), (d_+f -
d'_+f)\otimes d_-g>(x)
\\
&+&<(Ad^*_b\otimes Ad_b)(r_D), (d_-f -
d'_-f)\otimes d_-g>(x)
\\
&=& <(Ad^*_b\otimes Ad_b)(r_D), (d_+f -
d'_+f)\otimes d_-g>(x)
\end{eqnarray}
Here $d_+f$ is the differential "in the $G$-direction in $D(G)$",
$d_-f$ is the differential "in $G^*$-direction in $D(G)$" and
$r_D$ is the $r$-matrix for the double $D(\g)$ which is invariant
with respect to $Ad^*\otimes Ad$ action of $G^*$. Thus, the
Poisson bracket of two $G^*$-invariant functions is again
$G^*$-invariant and this proves the theorem.

Define the variety $D(G)//Ad_{D(G)}$ again as the geometric
quotient, i. e. as the spectrum of the $Ad_{D(G)}$-invariant
functions on $D(G)$. Since $D(G)=G\times G$ and $G//Ad_G\simeq
H/W$ we have
\[
D(G)//Ad_{D(G)}\simeq H/W\times H/W
\]
The natural imbedding $G^*\subset D(G)$ gives the map
\[
D(G)//{Ad_{G^*}}\to D(G)//{Ad_{D(G)}}
\]
and therefore the inclusion $C_{D(G)}(G)\subset C_{G^*}(G)$.

\begin{proposition}  The map $D(G)//{Ad_{G^*}}\to D(G)//{Ad_{D(G)}}$ is Poisson.

\end{proposition}

{\bf Proof}. We should show that the subalgebra
$C_{D(G)}(D(G))\subset C_{G^*}(D(G))$  is in the center of the
Poisson algebra $C_{G^*}(G)$. Let $f$ be an $Ad_{D(G)}$-invariant
function on $D(G)$ and $g$ be $Ad_{G^*}$-invariant function on
$D(G)$. Then $df=d'f, \ d_-g=d'_-g$ and therefore
\[
\{f,g\}(x)=<r_D, (d_+f-d'_+f)\otimes d_-g>=0
\]
Thus, the pull-back of $\eta$ gives central functions on
$D(G)/{Ad_{G^*}}$. This proves the proposition.

Consider the composition map

$\psi$ :

\begin{equation}\label{psi}
\psi : G\to D(G)\to D(G)//Ad_{G^*}
\end{equation}

The map $\psi$ is a composition of Poisson maps and therefore is a
Poisson map itself.
Now assume that $G$ is a factorizable Poisson Lie group.
In this case $D(G)=G\times G$.

\begin{theorem} \begin{enumerate}
\item $\psi^{-1}(\psi(g))=Z(g)\cap G'$ where $Z(g)$ is the
centralizer of $g$ in $G$ and $G'$ is the subset of factorizable
elements.

\item  The following diagram is commutative
\[
\begin{array}
{ccc} {D(G)}//Ad_{G^*} & \rightarrow & {D(G)}//Ad_{D(G)}=G//Ad_{G}
\times G//Ad_{G} \\ \uparrow        & {} & \uparrow \\ G &
\rightarrow & G//Ad_{G}
\end{array}
\]
Here left vertical arrow is the map $\psi$ and the right vertical
arrow is the diagonal embedding.

\end{enumerate}

\end{theorem}
Proof.

Now assume that $G$ is a simple Lie group with a factorizable
Poisson Lie structure. In this case $G//{Ad_G}\simeq H/W$ and we
have a composition of Poisson projections
\begin{equation}\label{G-psi-pi}
G\to \psi(G) \to H/W
\end{equation}

\begin{corollary} If $g$ simple element
$dim(\psi^{-1}(\psi(g)))=r$.
\end{corollary}

\subsubsection{The map $\beta$}

The following lemma is a combination of well known facts(see
\cite{Drinfeld87} and \cite{STS} for example).

\begin{lemma} \label{action} \begin{enumerate}
\item The subset $[B^-\times B^+]=\{(b_-,b_+)\in B^{-}\times
B^+|[b_+]_0=[b_-]_0\} $ is a Poisson submanifold in $B^-\times
B^+$.

\item The map
\begin{equation} \label{beta}
[B^-\times B^+]\to G, \ (b^-,b^+)\to b^+(b^-)
\end{equation}
is Poisson. Its image is open dense in $G$ and it is a covering
map with the group of deck transformations isomorphic to
$\Gamma=\{\varepsilon\in H| \varepsilon^2=1 \}$.

\item The adjoint action $h:(b_-,b_+)\mapsto (h^{-1}b_-h,hb_+h^{-1})$
of $j(H)$ on $[B^-\times B^+]$ is Poisson.
\end{enumerate}
\end{lemma}

Notice that $[B^-\times B^+]$ can be naturally identified with the
left coset $(B^-\times B^+)/j(H)$ as a Poisson manifolds.

\begin{corollary}
The coset manifold $[B^-\times B^+]/{Ad_{j(H)}}$  carries natural
Poisson structure.
\end{corollary}

\begin{proposition}
There is an isomorphism of Poisson manifolds
\begin{equation}\label{piso}
[B^-\times B^+]/Ad_{j(H)}\simeq [N^-\times N^+]/Ad_{j(H)}\times H
\end{equation}
\end{proposition}

This proposition follows directly from the explicit form of the
Poisson brackets for $B^+$ (see Appendix 2).

We have a natural map $B^-\times B^+\to G/Ad_{N^+}\times
G/Ad_{N^-}$ acting as
\begin{equation} \label{betas}
(b_-,b_+) \mapsto (Ad_{N^+}(b_-), Ad_{N^-} (b_+))
\end{equation}

This map induces the map
\begin{equation}\label{betaa}
\beta: [B^-\times B^+]/Ad_{j(H)}\to
D(G)//Ad_{G^*}=(G//Ad_{N^+}\times G//Ad_{N^-})/Ad_{j(H)}.
\end{equation}

\begin{theorem}

The map $\beta$
is  Poisson.
\end{theorem}
Proof. It follows from the theorem \ref{quot} that the map
$G\times G \to (G\times G)//Ad_{G^*}$ is Poisson. Consider the
composition map $ B^-\times B^+ \to G\times G \to (G\times
G)//Ad_{G^*}$. Each of these maps is Poisson, therefore the
composition is also Poisson.


In Appendix 2 we explain why maps $B^{\pm}\to G/Ad_{N^{\mp}}$ are
branched cover maps over their images. The third part of the
theorem follows from this.

Let $I: [B^-\times B^+]/{Ad_H}\to H$ be the projection to the
second factor in (\ref{piso}). It is a Poisson map.

\begin{theorem}\label{10}
\begin{enumerate}
\item Let $g\in G$ be a semisimple element (i.e. the elements which is
conjugate to an element of $H$). Then $\psi (g)\subset
\beta([B^-\times B^+]/Ad_{j(H)})$ (in other words, then there
exist elements $n_{\pm}\in N^{\pm}$ and $b_{\pm}\in B^{\pm}$ such
that $g=n_+b_-n_+^{-1}=n_-b_+n_-^{-1}$) .

\item The following diagram is commutative

\[
\begin{array}{ccc}
D(G)//Ad_{G^*} & \longrightarrow& D(G)//Ad_{D(G)}=H/W \times H/W
\\ \uparrow\be & {} & \uparrow \\ {[B^-\times B^+]}/Ad_{j(H)} &
\stackrel{\hat{I}}{\longrightarrow} & H/W
\end{array}
\]
\end{enumerate}
Here the map $\hat{I}$ is the composition of the map $I$ and of
the natural projection $H\to H/W$ and the right vertical map is
the the diagonal embedding $H/W\to H/W\times H/W$.
\end{theorem}

Proof. The second part of the theorem is obvious.

Let us prove that $\psi$ map semisimple elements into the image of
$\beta$ in $D(G)/Ad_{G^*}$. Assume that $x=ghg^{-1}$ where $h\in
H$ and $g\in G$. Assume that $g\in B^+wB^+$ for some $w\in W$. Let
$\bar{w}$ be a representative of $w$ in $N(H)\subset G$. Then
$g=b_+\bar{w}n_+ $ for some $b_+\in B^+$ and $n_+\in N_w^+$ where
$N_w^+$ is the subset of elements in $N^+$ which map into $N^-$
after conjugation with $\bar{w}$. We have:
\begin{eqnarray}
ghg^{-1}&=&b_+\bar{w}n_+hn_+^{-1}\bar{w}^{-1}b_+^{-1} \\
 &=& b_+ w(h) (\bar{w}h^{-1}n_+hn_+\bar{w}^{-1})b_+^{-1}  \\
&=& b_+w(h)\tilde{n}_-b_+^{-1}
\end{eqnarray}

This proves that each semisimple element can be written as
$x=n_+b_-n_+^{-1}$ for some $n_+\in N^+$ and $b_-\in B^-$.
Similarly one can prove that $x=n_-b_+n_-^{-1}$ for some $n_-\in
N^-$ and $b_+\in B^+$. This proves the first part of theorem.




\subsection{Symplectic leaves of ${\Bx}/{Ad_{j(H)}}$}

Let $u\in W$ be an element of the Weyl group and $u=s_{i_1}\dots
s_{i_\ell}$ be its reduced decomposition. The subset $|u|\subset
\{1,\dots ,r\}$ of numbers which appear in the sequence
$\{i_1,\dots ,i_r\}$ is called the support of $u$.

Denote by $H(u)$ the subgroup of $H$ generated by 1-parametric
subgroups corresponding to simple roots $\a_i$ with $i\in |u|$. We
have the following decomposition of $B^-\times B^+$:
\[
B^-\times B^+=\sqcup_{u,v \in W} B^-_u\times B^+_v ,
\]
where $B^{\pm}_w=B^{\pm}\cap B^{\mp}wB^{\mp}$. This decomposition
gives the decomposition of $[B^-\times B^+]$:
\[
[B^-\times B^+]=\sqcup_{u,v \in W} [B^-\times B^+]_{v,u}
\]
where
\begin{equation*}
{\Bx}_{v,u} =\{(b^-,b^+)\in B^-_v\x B^+_u\mid [b^+]_0 = [b^-]_0 ,
[b^{\pm}]_0\in H\}
\end{equation*}

\begin{theorem}\begin{enumerate}
\item The subsets ${\Bx}_{v,u}$ are Poisson submanifolds.

\item The Poisson submanifold ${\Bx}_{v,u}$ is fibered over
${\C}^{d_{u,v}}$ where $d_{u,v}=dim(ker_{{\mathfrak h}^*}(uv^{-1}-id)$
with fibers being common
level sets of functions $c_{u,v,t}(b^-b^+)$,
where $t\in ker_{\Lambda}(uv^{-1}-1)\subset {\mathfrak h}^*$.
Symplectic leaves of ${\Bx}_{v,u}$ are irreducible components of
these fibers.
\end{enumerate}
\end{theorem}

{\bf Proof.} The map ${\Bx}\to G$, $(b^-,b^+)\mapsto b^-b^+$ is
Poisson. It is a cover map with the group of deck transformations
$\Gamma=\{\e\in H\mid \e^2=1\}$ with the image which is open dense
in $G$. Symplectic leaves of $G$ are irreducible components of
level surfaces of functions $c_{u,v,t}$. The intersection of each
symplectic leaf of $G$ with the image of this map is open dense in
the symplectic leaf . Therefore symplectic leaves of ${\Bx}$ are
irreducible components of preimages of open dense subsets of
symplectic leaves of $G$. This proves the theorem.

\begin{proposition} The subgroup $(H(u) H(v))^{\perp}\subset H$ is the
stabilizer of the adjoint action of $j(H)$ on ${\Bx}_{v,u}$.
\end{proposition}

This follows from the factorization formulae for $B^-_v,B^+_u$
(see \cite{FZ}).

\begin{proposition} The  adjoint action of $H(u,v)=H/(H(u)
H(v))^{\perp}\simeq H(u) H(v)$ is transitive on ${\Bx}_{v,u}$ and
is Hamiltonian. Corresponding Hamiltonian vector fields are
generated by linear functions on the Lie algebra of $H(u,v)$.
\end{proposition}

Transitivity of the adjoint action of $H(u,v)$ is obvious. The
second part of the theorem follows from the $r$-matrix form of
Poisson brackets on $G$.

For $t=\sum_{i=1}^r t_i\omega_i\in ker_{\Lambda}(uv^{-1}-id)\subset
{\mathfrak{h}^*}$ define functions $[c_u\otimes c_v]_t$ on
$[B^-\times B^+]$ as
\begin{equation} \label{coc}
[c_v\otimes
c_u]_t(b_-,b_+)=\prod_{i=1}^r\D_{\omega_i,u^{-1}\omega_i}(b_+)^{t_i}
\D_{v\omega_i,\omega_i}(b_-)^{u^{-1}(t)i}
\end{equation}
These functions are $Ad_{j(H)}$-invariant and therefore define
functions on $[B^-\times B^+]/Ad_{j(H)}$.

\begin{corollary} The set of orbits of adjoint action of $j(H)$ passing
through ${\Bx}_{v,u}$ is isomorphic to
${\Bx}_{v,u}/Ad_{j(H(v,u))}$.
\end{corollary}

\begin{theorem}\begin{enumerate}
\item The isomorphism (\ref{piso}) induces the isomorphism of Poisson
varieties ${\Bx}_{v,u}/Ad_{j(H(u,v)}\simeq (N^-_v\times
N^+_u)/Ad_{j(H(v,u))}\times
H$. The Poisson structure on the last factor
is trivial.

\item Functions $[c_u\otimes c_v]_t$ are constant the subspace $H(u,v)$
in the second
factor.

\item Symplectic leaves of $(N^-_v\times N^+_u)/Ad_{j(H(v,u))}$ are
common level sets of functions $[c_u\otimes c_v]_t$, where $t\in
ker_{\Lambda}(uv^{-1}-1)\subset {\mathfrak h}^*$.
\end{enumerate}
\end{theorem}

This theorem follows from the Hamiltonian reduction via moment
map.

\subsection{Symplectic leaves of $\psi(G)\subset
D(G)//Ad_{j(G)}$}

Let $\psi: G\to (G\times G)/Ad_{j(G^*)}$ and $\beta: [B^-\times
B^+]/Ad_{j(H)}\to (G\times G)/Ad_{j(G^*)}$ be the maps defined in
(\ref{psi})(\ref{betaa}).

For $u,v\in W$ and $t=\sum_{i=1}^rt_i\omega_i\in {\mathfrak{h}^*}$
define functions $\tilde{c}_{u,v,t}$ on $G\times G$ as

\[
\tilde{c}_{u,v,t}(g_1,g_2)=\prod_{i=1}^r\D_{\O_i,u^{-1}\O_i}(g_1)^{t_i}
\D_{v\O_i,\O_i}(g_2)^{u^{-1}(t)_i}
\]

Here functions $\D_{u\O_i,v\O_i}$ are defined in (\ref{DDelta}).

\begin{proposition}\label{pull-b}
\begin{enumerate}
\item Function $\tilde{c}_{u,v,t}$ are invariant with respect to
the $Ad_{j(G^*)}$-action.
\item $c_{u,v,t}=\psi^*(\tilde{c}_{u,v,t})$.
\item $[c_u\otimes c_v]_t=\beta^*(\tilde{c}_{u,v,t})$ were
functions $[c_u\otimes c_v]_t$ are defined in (\ref{coc}).
\end{enumerate}
\end{proposition}
Proof. Recall that functions $\D$ have the following property:
\[
\D_{\O_i,u^{-1}\O_i}(\xi_-^{-1}g\xi_-)=[\xi_-]_0^{-\O_i+u^{-1}\O_i}
\D_{\O_i,u^{-1}\O_i}(g)
\]
\[
\D_{v\O_i,\O_i}(\xi_+^{-1}g\xi_+)=[\xi_+]_0^{\O_i-v\O_i}\D_{v\O_i,\O_i}(g)
\]
This proves the first part of the proposition. The second part
follows from the definition of $c_{u,v,t}$. For
$\beta^*(\tilde{c}_{u,v,t})$ we have:
\[
\beta^*(\tilde{c}_{u,v,t})(b_-,b_+)=\prod_{i=1}^r\
\D_{\O_i,u^{-1}\O_i}(b_+)^{t_i}
\D_{v\O_i,\O_i}(b_-)^{u^{-1}(t)_i}= [c_u\otimes c_v]_t(b_-,b_+)
\]
This proves the third part.

\bs Let $S_g^{u,v}\subset G$ be the symplectic leaf in $G$ passing
through $g$. Denote by $H^W(u,v)$ the set of $W$-orbits in $H$
which intersect $H(u,v)\subset H$.

\begin{theorem} Let $g=n_+b_-n_+^{-1}=n_-b_+n_-^{-1}\in G^{u,v}$ be
a semisimple element. Denote $[g]=Ad_{j(H)}(b_-,b_+)\in
[B_u^-\times B_v^+]/Ad_{j(H)}$ and denote by $\Sigma_{[g]}^{u,v}$
the symplectic leaf in $[B_u^-\times B_v^+]/Ad_{j(H)}$ passing
through $[g]$ . Let ${\mathcal O}_g$ be the $Ad_G$-orbit passing
through $G$. Then

\begin{enumerate}
\item if $(b_-',b_+')$ is a different pair representing $g$ then
$\Sigma_{[g]}^{u,v}=\Sigma_{[g]'}^{u,v}$.
\item $\psi(S_g^{u,v}\cap {\mathcal O}_g)=\beta(\Sigma_{[g]}^{u,v})$.
\item $\pi\circ\psi (S^{u,v}_g\cap {\mathcal O}_g)=A_G(g)\in H^W(u,v)$,
assuming the identification $G//Ad_G=H/W$.
\end{enumerate}
\end{theorem}
Proof. The symplectic leaf $S^{u,v}_g$ is the common level surface
of functions $c_{u,v,t}$ which contains $g$. It follows from this
fact and from the proposition \ref{pull-b} that the common level
surface of function $\tilde{c}_{u,v,t}$ which contains $\psi(g)$
is a Poisson subvariety in $\psi(G)$. On the other hand
$\psi(g)=\beta([g])$ and by similar reasoning the subvariety of
the common level surface of $\tilde{c}$ passing through $\psi(g)$
which consists of $\psi$ images of elements
$g'=n'_+b'_-{n'}_+^{-1}=n'_-b'_+{n'}_-^{-1}$ with
$[b'_+]_0=w[b_+]_0$ for some $w\in W$, is
$\beta(\Sigma^{u,v}_{[g]})$.

\begin{corollary} If $g$ is semisimple the submanifold
$\psi(S_g^{u,v}\cap {\mathcal O}_g)=\beta(\Sigma^{u,v}_{[g]})$ is a symplectic leaf
in $\psi(G)$ .
\end{corollary}
\begin{corollary}\label{2} If $[h]\in H^W(u,v)$ is a generic orbit
corresponding to the coset $Ad_G(g)$,
then connected components of $\pi^{-1}([h])$ are symplectic leaves of
$\psi(S^{u,v}_g)$ of dimension $dim(S^{u,v})-2dim(H(u,v)$.
\end{corollary}

\subsection{Integrability of characteristic systems on $G$}

In order to prove the integrability of the characteristic system
on the symplectic leaf $S^{u,v}\subset G^{u,v}$ we should describe the system of projections
with properites (\ref{psi-pi}).

Let $S^{u,v}_g$ be a symplectic leaf in $G$ through $g\in G^{u,v}$.
The restriction of Poisson maps
(\ref{G-psi-pi}) to $S^{u,v}_g\subset G$
 gives the composition map
\begin{equation}\label{S-psi-pi}
S^{u,v}_g\stk{\psi}{\longrightarrow}
\psi(S^{u,v}_g)\stk{\pi}{\longrightarrow} Ad_G(S^{u,v}_g)\simeq H^W(u,v)
\end{equation}
where $Ad_G(S^{u,v}_g)$ is the set $Ad_G$ orbits intersecting
$S_{u,v}$ and $H^W(u,v)$ is the set of $W$-orbits in $H$ passing through
$H(u,v)$.

Let $[h_g]\in H^W(u,v)$ be the element corresponding to the orbit
$Ad_G(g)$. According to the corollary \ref{2} connected
components
of $\pi^{-1}([h_g])$, $[h_g]\in H^W(u,v)$ are symplectic leaves of
$\psi(S^{u,v}_g)$ of
dimension
\[
\dim(\pi^{-1}([h_g]))=\dim(S^{u,v}_g)-2\dim(H(u,v)) \ .
\]

\begin{lemma} Let $g\in G^{u,v}$ be a semisimple element.
Then $dim(\psi^{-1}(\psi(g))=dim(H(u,v)$.
\end{lemma}

Proof. Assume that $g.g'\in G^{u,v}$ are semisimple elements such that
$\psi(g)=\psi(g')$ Then
\[
g'=\tilde{b_\pm} g\tilde{b_\pm}^{-1}
\]
for some $\tilde{b_\pm}\in B^{\pm}$ with $[\tilde{b_+}]_0=[\tilde{b_\pm}]_0^{-1}$.

Since $g,g'$ are semisimple we can represent them as
\[
g=n_{\pm}b_{\mp}n_{\pm}^{-1}, \ \ g'=n'_{\pm}b'_{\mp}{n'}_{\pm}^{-1}
\]
where $b_-\in B^-_v$ and $b_+\in B^+_u$.
Then we have:
\[
b'_-=\beta_+b_-\beta_+^{-1}, \ \ \beta_+={n'}_+^{-1}\tilde{b}_+n_+,
\]
\[
b'_+=\beta_-b_+\beta_-^{-1}, \ \ \beta_-={n'}_-^{-1}\tilde{b}_-n_-.
\]
Notice that $[\beta_\pm]_0=[\tilde{b}_\pm]_0$ and therefore
$[\beta_+]_0=[\beta_-]_0^{-1}$. Let $\beta_\pm=h^{\pm 1}\nu_\pm$
where $\nu_\pm\in N^\pm$. For $\nu_\pm$ there is only discrete
choice (determined by the action of the Weyl group on $H$). The
subgroups $H(u)^{\vee}$ and $H(v)^{\vee}$ act trivially ( via
conjugation) on $B^+_u$ and on $B^-_v$ respectively. Therefore for
given semisimple $g\in G^{u,v}$ the variety of semisimple elements
$g'\in G^{u,v}$ such that $\psi(g)=\psi(g')$ has dimension
$dim(H(u,v)$ and therefore $dim(\psi^{-1}(\psi(g)))=dim (H(u,v)$.

This lemma together with previous results proves the following
theorem.

\begin{theorem} Projection $\psi$  in (\ref{S-psi-pi}) has
$dim(H(u,v)$-dimensional kernel, the image of $\pi$ has the same
dimension and connected components of $\pi^{-1}$ of generic points
are symplectic leaves in the image of $\psi$. Therefore a
Hamiltonian system generated by an $Ad_G$-invariant function on
$S^{u,v}_g$ is integrable.
\end{theorem}

\begin{rem} When $H(u,v)\neq H$ (or, equivalently, when reduced
decompositions of $u$ and $v$ contain all simple reflections) we
say that the symplectic leaf $S^{u,v}$ is not of full rank. In
this case it is a symplectic leaf of the full rank an appropriate
semi-simple Poisson Lie subgroup in $G$.
\end{rem}
\begin{rem}
Among symplectic leaves in $G$ of full rank there are  symplectic
leaves corresponding Coxeter elements. They have dimension $2r$
and are integrable in the usual Liouville sense (the invariant
tori have dimension $r$). Corresponding integrable systems have
been studied in \cite{HKKR}. They are deformations of Toda
systems.
\end{rem}

\subsection{Action-angle variables}

We will say the element $x \in G$ of a simple, complex Lie group
is generic if it is a conjugate to a generic element from the Cartan subgroup:
$x = u h u^{-1}.$ Let $V_{\lambda}$ be a finite dimensional irreducible
representation of $G$ with the weight decomposition
$$V_{\lambda} = \oplus_{\mu \in \Di(\lambda)} V_{\lambda}(\mu).$$

For generic $x \in G$ denote by $P_{\mu}^{\lambda}$ the complete
system of orthogonal projections on the eigenspace of $x$ in
$V_{\lambda}:$
$$
x = \sum_{\mu} t_{\mu} P_{\mu}^{\lambda}, \qquad
 P_{\mu}^{\lambda} P_{\nu}^{\lambda} =
P_{\mu}^{\lambda} \delta_{\mu,\nu}
$$
where $t_{\mu} \in C^{\times}.$ Since $x$ is generic
$P_{\mu}^{\lambda} = u Q_{\mu}^{\lambda} u^{-1},$ where
$Q_{\mu}^{\lambda}$ is the projection to the subspace of
$V_{\lambda}(\mu)$ in the weight decomposition of $V_{\lambda}.$
For the same reason $t_{\mu}$ is the value of $h $ on $V_{\lambda}(\mu).$

Let $H$ be an $Ad$-invariant function on $G$ and $g_{\pm}(t,x)$
be the factorized components of $g(t,x) = \exp{t \nabla H(x)}:$
\begin{equation} \label{decomp}
g_+(t,x) g_-(t,x)^{-1} = g(t,x)
\end{equation}
This factorization exists for sufficiently small $t.$

Denote by $v_{\lambda}$ the highest weight vector of $V_{\lambda}$
and by $(\cdot,\cdot)$ the Cartan form, i.e the non degenerate
bilinear form such that $(\omega(a) x, y) = (x, a y)$ where
$\omega$ is the Cartan antiinvolution
($[\omega(a),\omega(b)]=-\omega([a,b])$):
\[
\omega(e_i)=f_i, \  \omega(f_i)=e_i, \ \omega(h_i)=h_i
\]
and we assume the normalization $(v_{\lambda}, v_{\lambda})=1.$

Introduce variables (functions on generic elements of $G$):

$$ r_{\mu}^{\lambda} = (v_{\lambda}, P_{\mu}^{\lambda}
v_{\lambda})
$$

\begin{theorem}
Let $\{x(t)\}$ be the flow line of the Hamiltonian vector field
generated by H, passing through c at t=0. Then
\begin{equation}
\label{reshi1}
r_{\mu}^{\lambda}(x(t)) = \frac{e^{-t X_{\mu}(x)} r_{\mu}^{\lambda}}{
\sum_{\nu \in \Di_{\lambda}} e^{-t X_{\nu}(x)} r_{\nu}^{\lambda}}
\end{equation}
Here $X_{\mu}(x)$ is the eigenvalue of $\nabla H(x)$ on $P_{\mu}^{\lambda}.$
\end{theorem}

\noindent {\em Proof:} According to \cite{STS} we have: $$ x(t) =
g_{+}(t,x)^{-1} x g_{+}(t,x) $$ where $g_{\pm}(x,t)$ are
determined by (\ref{decomp}). Therefore $$ r_{\mu}^{\lambda}(x(t))
= (v_{\lambda}, P_{\mu}^{\lambda}(x(t))v_{\lambda}) =
(v_{\lambda},
g_{+}(t,x)^{-1}P_{\mu}^{\lambda}(x(t))g_{+}(t,x)v_{\lambda}) $$ On
the other hand $g_{\pm}(t,x)$ are elements of the Borel
subalgebras $B_{\pm}$ of $G$. Write $g_{\pm}(t,x)$ as
\begin{eqnarray*}
g_{+}(t,x) &=& u_{+}(t,x) h(t,x) \\
g_{-}(t,x) &=& u_{-}(t,x) h(t,x)^{-1}
\end{eqnarray*}
where $u_{\pm}$ belong to corresponding unipotent subgroups and
$h$ is in the Cartan subgroup.

According to the definition of $(\cdot,\cdot)$ we have: $$
r_{\mu}^{\lambda}(x(t)) = (\omega(g_{-}(t))^{-1} v_{\lambda},
g(t)^{-1} P_{\mu}^{\lambda} g_{+}(t) v_{\lambda}) $$ and therefore
\begin{equation}
\label{reshi1a}
r_{\mu}^{\lambda}(x,t) = h_{\lambda}(t) \bar{h}_{\lambda}(t)
e^{t X_{\mu}(x)} r_{\mu}^{\lambda}
\end{equation}
where
$$ h(t) v_{\lambda} = h_{\lambda}(t) v_{\lambda}, \qquad
\omega(h(t)) v_{\lambda} = \bar{h}_{\lambda}(t)v_\lambda
$$
On the other hand
\begin{equation}
\label{reshi1b}
( v_{\lambda}, g(t)^{-1} v_{\lambda}) = (v_{\lambda}, g_{-}(t)
g_+(t)^{-1} v_{\lambda}) = h_{\lambda}(t)^{-1} \bar{h}_{\lambda}(t))^{-1}
\end{equation}
and
$$
( v_{\lambda}, g(t)^{-1} v_{\lambda}) = \sum_{\mu \in \Di_{\lambda}}
e^{-X_{\mu}(x)t} r_{\mu}^{\lambda}
$$
substituting this into (\ref{reshi1a}) and (\ref{reshi1b}) we obtain
(\ref{reshi1}).
\eb

Let $\mu_1$ and $\mu_2$ be two weights in an irreducible
representation $V_{\lambda}.$ Consider $$ r_{\mu_1,\mu_2} =
\frac{r_{\mu_1}^{\lambda}}{r_{\mu_2}^{\lambda}} $$ We have: $$
r_{\mu_1 \mu_2}(t) = exp(t (X_{\mu_2}(x)- X_{\mu_1}(x))) r_{\mu_1
\mu_2} $$ Therefore the logarithms of $r_{\mu_1 \mu_2}$ are affine
coordinates on invariant tori and therefore $n$ independent
variables of this type can serve as angle variables. The
eigenvalues of $x$ are action variables for the Toda system.

For example for $SL_n$ we can choose $\lambda=\omega_1,
\mu=\omega_1-\alpha_1-\dots-\alpha_i$. This is equivalent to
Moser's construction for Toda symplectic leaf for $SL_n$.

\begin{rem} Coxeter-Toda systems are characteristc systems on symplectic
leaves corresponding to a pair of Coxeter elements of the Weyl
group. In real totally positive case the action-angle variables
are global coordinates on the phase space for such systems. It
will be interesting to see if similar property holds for any
symplectic leaf.
\end{rem}

One should note that this construction of action-angle variables
is very similar to the one given by Kostant \cite{Kostant-oxford}
in linear case.

\section{Conclusion}

We proved that a Hamiltonian system on any symplectic leaf of a
simple Poisson Lie group with the standard Poisson structure is
integrable if the Hamiltonian is $Ad_G$-invariant. Liouville tori
of such systems are intersections of dressing and adjoint orbits.

One of the most interesting next questions is to describe the
spectrum of corresponding quantum systems. In case of Toda systems
this involves Wittaker vectors and some other facts about
principal unitary series of representations of split real form of
$G$.

In our case the analongs of Wittaker vectors and of principal
unitary series of representations for the split real forms of
quantized universal enveloping algberas $U(\mathfrak g)$ should
play similar role.

\section{Appendix 1.}
Here will prove some useful fact which was not used this paper.

\begin{theorem} Let $(D(G),p)$ be the double of Poisson Lie group $G$.
Let $(D(G),p_*)$ be the Poisson structure on the manifold $D(G)$
induced by the factorization map. Then if $f$ and $g$ are
$Ad_{G^*}$ invariant functions on $D$,
\[
\{f,g\}=\{f,g\}_* \  .
\]
\end{theorem}

{\bf Proof.} The Poisson brackets $\{\dt,\dt\}$ and  $\{\dt,\dt\}_*$
have the following form
\begin{eqnarray*}
\{f,g\} &=& \<r,df\o dg\>-\<r,d'f\o d'g\>
=(d_+f, \ d_-g)-(d'_+f, \ d'_-g) \\
\{f,g\}_* &=& +\<r,df\wedge dg+d'f\wedge d'g\> +
\<r,d'f\o dg-d'g\o df\> \\
&=& +\tfrac 12((d_+f, \ d_-g)+(d'_+f, \ d'_-g)-(d_+g, \ d_-f)-
(d'_+g, \ d'_-f)) \\
&&\quad -(d'_+f, \ d_-g)+(d'_+g, \ d_-f)
\end{eqnarray*}
Now assume that $f$ and $g$ are $Ad_{G^*}$-invariant. Then
$d_-f-d'_-f$, \ $d_-g=d'_-g$ and we have
\begin{eqnarray*}
\{f,g\} &=&
(d_+f-d'_+f, \ d_-g)=\tfrac 12(d_+f-d'_+f, d_-g)-
\tfrac 12(d_+g-d'_+g, \ d_-f)\\
\{f,g\}_* &=&\tfrac 12((d_+f,d_-g)-(d'_+f,d_-g)-(d_+g,d_-f)+
(d_+g,d'_-f) \\
&=&-\tfrac 12 (d_+f,d'_-g-d_-g)+\tfrac 12
(d_+g,d'_-f-d_-f)
\end{eqnarray*}
The theorem follows.

\section{Appendix 2. The Poisson structure on $G//Ad_{B^-}$}

\subsection{Poisson brackets of $Ad_{B^-}$-invariant functions on
$G$} Recall that the Poisson bracket on functions on the double
$D(B_+)=G\x H$ of $B_+$ (with standard Poisson structure) has the
following form
\[
\{f,g\} = \<r_0,d^+_0f\o d^-_0g\> -\<r_0,d^{+'}_0f\o d^{-'}_0g\>
+\<r_1,d_+f\o d_-g\> - \<r_q,d'_+f\o d'_-g\> \ .
\]
Here we assume that decompositions $B_\pm=HN_{\pm}$ are fixed,
together with imbeddings $B_{\pm}\hookrightarrow G\x H$, \
$b_{\pm}\mapsto (b_{\pm},([b]_0)^{\pm 1})$, where $[ \
]_0:B_{\pm}\mapsto H$ are projections to the Cartan subgroup.

Elements $r_0,r_1$ are canonical elements in ${\f}\o{\f}$ and in
${\N}_+\o{\N}_-$ respectively
\[
r_0=\sum_i H^+_i\o (H^-)^i \ , \qquad r_1=\sum_{\a}e_{\a}\o f_{\a}
\ .
\]
Here we assume that the first factor in ${\f}\o{\f}$ is the image
of the Cartan subalgebra in $B$  and the second in the image of
the Cartan subalgebra in $B_-$ under projections $[ \ ]_0$.

Differentials $d^{\pm}_0$ are taken in the direction of Cartan
subgroup $H\subset B^{\pm}$ imbedded in $G\x H$ via maps $i$ and
$j$ respectively . Differentials $d_{\pm}$ are taken in the
direction of $N_{\pm}\subset G\x H$. If one trivializes tangent
bundle to $G$ by identifying tangent spaces to $G$ with $\G$ we
have
\begin{eqnarray*}
\<\xi_0,d^+_0f\>(g,h) &=&\frac{d}{dt} f(e^{t\xi_0}g, \
e^{t\xi_0}h)|_{t=0} \ , \\ \<\xi_0,d_0^{+'}f\>(g,h)
&=&\frac{d}{dt} f(ge^{t\xi_0}, \ he^{t\xi_0})|_{t=0} \ , \\
\<\xi_"0,d^-_0f\>(g,h) &=& \frac{d}{dt} f(e^{t\xi_0}g, \
e^{t\xi_0}h)|_{t=0} \ , \\ \<\xi,d^{-'}_0f\> (g,h) &=&\frac{d}{dt}
f(ge^{t\xi_0}, \ he^{-t\xi_0})|_{t=0}
\end{eqnarray*}
and
\begin{eqnarray*}
\<\xi_{\pm},d_{\pm}f\>(g,h) &=& \frac{d}{dt} f(e^{\xi_{\pm}t}g, \
h)|_{t=0} \ , \\ \<\xi_{\pm},d'_{\pm}f\>(g,h) &=& \frac{d}{dt}
f(ge^{\xi_{\pm}t}, \ h)|_{t=0} \ .
\end{eqnarray*}
Here $\xi_0\in{\f}$, \ $\xi_{\pm}\in{\N}_{\pm}$.

In other words,
\begin{eqnarray*}
\<\xi,d^{\pm}_0f\> &=& \<\xi,d_0f\>\pm \<\xi,d_Hf\>\\
\<\xi,d^{\pm}_0f\> &=& \<\xi,d'_0f\>\pm \<\xi,d_Hf\>
\end{eqnarray*}
where $d_0 f,d'_0f$ are left and right differentials of $f$ in the
Cartan direction $H\subset G$, \ $d_Hf$ is the differential in the
direction of the second factor in $G\x H$. Thus, for the Poisson
bracket we have
\begin{eqnarray*}
\{f,g\} &=& \<r_0,d_0f\o d_0g\>-\<r_0,d'_0f\o d_0g\> \\ &&\quad +
\<r_0,d_Hf\o d_0g-d_0f\o d_Hg\>- \<r_0,d_Hf\o d'_0g-d'_0f\o d_Hg\>
\\ &&\quad +\<r_1,d_+f\o d_-g\> - \<r_1,d'_+f\o d'_-g\>
\end{eqnarray*}
Antisymmetrizing this bracket we obtain:
\begin{eqnarray*}
\{f,g\} &=&
 \<r_0,d_Hf\o (d_0g-d'_0g)-(d_0f-d'_0f)\o d_Hg\>\\
&&\quad +\tfrac 12\<r_1,d_+f\o d_-g - d_+g\o d_-f - d'_+f\o d'_-g
+d'_+g\o d'_-f\> \ .
\end{eqnarray*}
Let $C_{Ad_{B^-}}(D)$ be the algebra of $Ad_{B^-}$-invariant
functions on $D$.

As it follows from the previous subsection the algebra
$C_{Ad_{B^-}}(D)$ is finitely generated. Define the variety
$D//Ad_{B^-}= Spec (C_{Ad_{B^-}}(D))$.

\bs The adjoint action of $B^-$ on $D$ is trivial on $H$-component
of $D=G\x H$. Thus,
\begin{equation}\label{star}
D//Ad_{B^-} = G//Ad_{B^-} \x H
\end{equation}
as a variety.

\begin{lemma}
The formula (\ref{star}) describes the Poisson variety
$D(G)//Ad_{B^-}$ as the product of two Poisson varieties with the
trivial Poisson structure on the $H$-factor.
\end{lemma}

Indeed, if $f$ and $g$ are $Ad_{B^-}$-invariant functions on
$G\times H$ we have \ $d_0f=d'_0f$, \ $d_-f=d'_-f$ and the same
for $g$. Therefore for the Poisson bracket between $f$ and $g$ we
have
\begin{equation}\label{**}
\{f,g\} =\tfrac 12\<r_1,(d_+f-d'_+f)\o d_-g\>- \tfrac
12\<r_1,(d_+g-d'_+g)\o d_-f\>
\end{equation}
This means that functions constant along $G$ are central in the
Poisson algebra which proves the lemma.

The map $B^+\hookrightarrow D(B^+)\to D(B^+)//{Ad_{B^-}}$ is a
composition of Poisson maps and therefore is Poisson. Projecting
to the second factor in the (\ref{star}) we have the Poisson
projection:
\begin{equation}\label{m}
B^+ \to G//{Ad_{B^-}} .
\end{equation}

The projection $B^+\to B^+/Ad_H$ is Poisson. This follows from the
$Ad_H$-invariance of the standard Poisson structure on $G$. It is
also clear that the diagram
\[
\begin{array}{ccc}
B^+ & \longrightarrow & G//Ad_{B^-} \\ [0.125in] {} & \searrow &
\uparrow\pi \\ [0.125in] {} & {} &  B^+/Ad_H
\end{array}
\]
is commutative. Therefore the map
\[
B^+\to G//Ad_{B^-}
\]
is Poisson.

The image of this map is open dense in $G//Ad_{B^-}$ and the map
is a finite branched cover. The number of branches over generic
point is equal to $|W|$ and the Weyl group $W$ acts naturally on
the fibers.

\subsection{Poisson structure on $B^+/Ad_H$}
The Poisson structure on $B^+/Ad_H$ can be described
explicitly.

\begin{theorem}\label{adjoint} The Poisson bracket of two $Ad_H$-invariant
functions on $B^+$ has the following form:
\begin{equation} \label{b}
\{f,g\}(x) =\<({\rm{id}}\o{\rm{Ad}}_{x^{-1}}(r_1), d_+f\o
\p'_+g\>(x)
\end{equation}
Here $d_+f$ is the left differential of $f$ at the point $x\in
B^+$ and $\p'_+g$ is the right differential of $g$ projected on
${\mathfrak n}_+\subset T_xB^+$.
\end{theorem}

Proof. This theorem can be derived as a pull-back of the Poisson
structure on $G//Ad_{B^-}$ or from the restriction of the standard
Poisson structure on $G$ to $Ad_H$-invariant functions.

Let us first compute it as a pull-back.

Assume that $x \in G$ belongs to the image of (\ref{m}), i.e.
there exists an element $x_+\in B^+$ and $n_-\in N^-$ such that
$x=n_-x_+n_-^{-1}$.

\begin{lemma} The value of the Poisson bracket of two
${\rm{Ad}}_{B^-}$-functions on such element $x=n_-x_+n_-^{-1}\in
G$ has the following form
\begin{equation}\label{a}
\{f,g\}(x_+) =\{f,g\}(x_+)=\sum_i \<\rho^i,d_+f-d'_+f\> (x_+)
\<\rho_i,d_-f\>(x_+)
\end{equation}
\end{lemma}
Proof. For $\xi_+\in{\N}_+$ and $b_-\in B^-$ we have
\begin{eqnarray*}
&& \<\xi_+,d_+f-d'_+f\>(b_-xb_-^{-1}) =\frac{d}{dt}
f(e^{\xi_+t}b_-xb_-^{-1} \ e^{-\xi_+t})|_{t=0} \\ &&\quad =
\frac{d}{dt} f(e^{(\xi_+^{b_-^{-1}})_+t} x \
e^{-(\xi_+^{b_-^{-1}})_+t})|_{t=0}
\\ &&\quad =\<(\xi_+^{b_-^{-1}})_+, \ d_+f-d'_+f\>(x)
\end{eqnarray*}
Here $(\xi_+^{b_-^{-1}})_+\in{\N}_+$. On the other hand,
\[
\<\xi_-,d_-g\>(b_-xb_-^{-1}) =\frac{d}{dt}
g(e^{t\xi_-}b_-xb_-^{-1})|_{t=0} =\frac{d}{dt}
g(e^{\xi_-^{b_-^{-1}}}x)|_{t=0} =\<\xi_-^{b_-^{-1}},d_-g\>(x_+) \
.
\]
The element $r_1=\sum_i \rho^i\o \rho_i$ is
$({\mb{Ad}}^*\o{\mb{Ad}})_{N^-}$-invariant and is also invariant
with respect to the diagonal action of $H$. Therefore
\begin{eqnarray*}
\<r_1,(d_+f-d'_+f)\o d_-g\>(b_-xb_-^{-1}) &=& \sum_i
\<(((\rho^i)^{b_-^{-1}})_+, d_+f-d'_+f\>(x) \<\rho_i^{b_-^{-1}},
d_-f\>(x_+) \\ &=& \sum_i \<\rho^i,d_+f-d'_+f\>(x)
\<\rho_i,d_-f\>(x)
\end{eqnarray*}

This proves the lemma.

To prove the theorem we should verify that the Poisson bracket
(\ref{a}) between two $Ad_{B^+}$-invariant functions is given by
(\ref{b}).

For $\xi_+\in{\N}_+$ and $x_+\in B^+$ we have:
\[
\<\xi_+,d_+f-d'_+f\>(x_+) =
\<\xi_+-{\mb{Ad}}_{x_+}(\xi_+),d_+f\>(x_+) \ .
\]
For $\xi_-\in{\N}_-$ and $x_+\in B_+$ we have
\begin{eqnarray*}
\<\xi_-,d_-g\>(x_+) &=&\frac{d}{dt} g(e^{t\xi_-}x_+)|_{t=0}\\
&=&\frac{d}{dt} g(e^{t\eta_-}e^{ta_+}x_+e^{-t\eta_-})|_{t=0} \\
&=&\frac{d}{dt} g(e^{ta_+}x_+)|_{t=0} =\<a_+,\p_+g\>(x_+)
\end{eqnarray*}
Here $\eta_-\in{\N}_-$ and $a_+\in b_+$ satisfy the equation
\[
\xi_-=\eta_-+a_+-{\mb{Ad}}_{x_+}(\eta_-) \ .
\]
This equation gives the equation for $\eta_-$:
\begin{equation}\label{sttr}
\xi_-=\eta_- -({\mb{Ad}}_{x_+}(\eta_-))_-
\end{equation}
and for $a_+$ we have
\[
a_+=  ({\mb{Ad}}_{x_+}(\eta_-))_+
\]
Here $({\mb{Ad}}_{x_+}(\eta_-))_+$ is the ${\b}_+$ part of
${\mb{Ad}}_{x_+}(\eta_-)$. Thus,
\[
\<\xi_-,d_-g\>(x_+) =\< ({\mb{Ad}}_{x_+}(\eta_-))_+, \
\p_+g\>(x_+)
\]
where $\eta_-\in{\N}_-$ is the solution to (\ref{sttr}). Here
$\p_+f$ is the differential of $f$ "in the direction of $B^+
\subset G$.

For the value of the Poisson bracket of two Ad${}_{B^-}$-invariant
functions on the element $x=n_-xn_-^{-1} \in G$ we have ,
\begin{eqnarray*}
\{f,g\}(x)=\{f,g\}(x_+)&=&\<r_1-({\mb{Ad}}_{x_+}\o{\mb{id}})(r_1),
d_+f\o d_-g\>(x_+) \\ &=&
\<\rho^i-{\mb{Ad}}_{x_+}(\rho^i),d_+f\>(x_+)
\<{\mb{Ad}}_{x_+}(\sig_i),\p_+g\>(x_+)
\end{eqnarray*}
Here $\sig_i$ is the solution to
\[
\rho_i=\sig_i-({\mb{Ad}}_{x_+}(\sig_i))_- \ .
\]
Because $r_1$ is the canonical element in $\N_+\o\N^-$ (assuming
we fixed an iso\mor \\ $\N^-\cong (\N_+)^*$ by the choice of
Killing form) we have:
\[
{\mb{Ad}}_{x_+}(\rho^i)\o\rho_i =\rho^i\o
({\mb{Ad}}_{x^{-1}_+}(\rho^i))_-
\]
and therefore
\[
(\rho^i-{\mb{Ad}}_{x_+}(\rho^i))\o\sig_i=-\rho^i\o
({\mb{Ad}}_{x_+^{-1}}(\rho^i))_- \ .
\]
Thus,
\[
\{f,g\}(x_+) =-\<\rho^i,d_+f\>(x_+)
\<{\mb{Ad}}_{x_+}(({\mb{Ad}}_{x_+^{-1}}(\rho^i))_-),\p_+g\>(x_+) \
,
\]
On the other hand,
\begin{eqnarray*}
&& \<{\mb{Ad}}_{x_+}({\mb{Ad}}_{x_+^{-1}}(\rho^i)_-),\p_+g\>(x_+)
\\ &&\quad =\<\rho_i-
{\mb{Ad}}_{x_+}({\mb{Ad}}_{x_+^{-1}}(\rho^i))_+),\p_+g\>(x_+) \\
&&\quad =-\< {\mb{Ad}}_{x_+^{-1}}(\rho^i)_+,\p'_+g\>(x_+)
\end{eqnarray*}
The theorem is proved.

\bs

\subsection{The second proof of the theorem \ref{adjoint}}

Let $G//N^-$ be the categorical quotient for the right action of
$N^-$ on $G$. The following is well known.

\begin{theorem} The map $G\to G//N^-$ is Poisson.
\end{theorem}
Proof. Let $f$ and $g$ be two functions on $G$ invariant with
respect to the right action of $N^-$. We have $d'_-g=0$ and
\[
<r_0,d'_0f\otimes d'_0g>(xn_-)=\sum_{i=1}^r\frac{d^2}{dsdt}
f(xn_-e^{th_i}) g(xn_-e^{sh^i})=<r_0,d'_0f\otimes d'_0g>(x)
\]
For the Poisson brackets of two such functions we have:
\[
\{f,g\}=\< r_0,d_0f\o d_0g-d'_0f\o d'_0g\>
 +\tfrac 12\< r_1,d_+f\o d_-g\>(x_+)-
\tfrac 12\<r_1,d_+g\o d_-f\>
\]
Here the first and the second term are invariant with respect to
the left action of $N^-$ the invariance of the second term is
shown above and the forth term vanishes. Thus, the Poisson bracket
of two invariant functions again invariant. Therefore the map is
Poisson.

\begin{theorem}\begin{enumerate}
\item The map $\phi: B^+\hookrightarrow G \to G//N^-$ is Poisson
for the standard Poisson structure on $G$ and and is a local
isomorphism ( an isomorphism in a neighborhood of the identity in
$B^+$) .
\item Let $f$ and $g$ be two left
$N^-$-invariant functions on $G$, then the pull-back of $\phi$
gives the following Poisson brackets
\begin{eqnarray*}
\{f,g\}&=&\< r_0,d_0f\o d_0g-d'_0f\o d'_0g\> \\ && +\tfrac
12\<({\rm{id}}\o{\rm{Ad}}_{x_+^{-1}})(r_1),d_+f\o\p_+g\> -\tfrac
12\<({\rm{id}}\o{\rm{Ad}}_{x_+^{-1}})(r_1),d_+g\o\p_+f\> \\ &=& \<
r_0,d_0f\o d_0g-d'_0f\o
d'_0g\>+\<({\rm{id}}\o{\rm{Ad}}_{x_+^{-1}})(r_1),d_+f\o\p_+g\>
\end{eqnarray*}
This Poisson bracket is the standard Poisson structure on $B^+$.

\end{enumerate}
\end{theorem}

{\bf Proof.} The map $\phi$ is Poisson since the since the left
action of $N^-$ on $G$ is admissible (in a sense of \cite{STS}).
It is also clear that it is an isomorphism in a neighborhood of
$1$.

For the Poisson bracket on functions on $G$ we have:
\[
\{f,g\}=\< r_0,d_0f\o d_0g-d'_0f\o d'_0g\>
 +\tfrac 12\< r_1,d_+f\o d_-g\>(x_+)-
\tfrac 12\<r_1,d_+g\o d_-f\>
\]
For \ $\xi\in\N^-$ \ we have
\[
\<\xi,d_-f\>(x)=\frac{d}{dt} f(e^{t\xi}x_+)|_{t=0} = \frac{d}{dt}
f(x_+ \, e^{t(\xi^{x_+^{-1}})_+})|_{t=0}
\]
Here we used decomposition $\eta=\eta_++\eta_-$ where $\eta=\G$, \
$\eta_+\in\b_+$, \ $\eta_-\in\N^-$. Thus,
\begin{eqnarray*}
\{f,g\}&=&\< r_0,d_0f\o d_0g-d'_0f\o d'_0g\> \\ && +\tfrac 12
\<({\rm{id}}\o{\rm{Ad}}_{x_+^{-1}})(r_1),d_+f\o\p'_+g\> \\ &&
-\tfrac 12 \<({\rm{id}}\o{\rm{Ad}}_{x_+^{-1}})(r_1),d_+g\o\p'_+f\>
\ .
\end{eqnarray*}

\begin{corollary} If functions $f$ and $g$ are $Ad_H$-invariant,
\[
\{f,g\}=\<({\rm{id}}\o{\rm{Ad}}_{x_+^{-1}})(r_1),d_+f\o\p_+g\>
\]
This gives the second proof of the formula (\ref{b}).
\end{corollary}

\end{document}